\documentclass[12pt,reqno]{amsart}
\usepackage[all]{xy}
\usepackage{graphicx}
\usepackage{amssymb}
\usepackage{eepic}
\makeatletter
\newdimen\sectionruledimen
\def\makeparrule{%
 \def\par{%
  \endgraf\nobreak\vskip\lineskip\nointerlineskip
  \hbox to\hsize{\hskip\sectionbackskip\leaders\hrule height
  \sectionruledimen\hfil}%
  }%
 }%

\def\section{%
 \@startsection {section}{1}{\sectionbackskip}{-10pt plus -1.2ex minus -.2ex}%
  {.5pt}{\normalsize\bf\makeparrule}%
 }%
\newdimen\sectionbackskip
\sectionbackskip\z@
\newdimen\sectionruledimen
\sectionruledimen\z@

\def\punteada{\leaders\hbox{$\m@th \mkern1.5mu - \mkern1.5mu$}\hfill}
\makeatother

\newtheorem{theorem}{Theorem}
\newenvironment{theor}{\smallskip\begin{trivlist}
 \item[\hspace{\labelsep}{\noindent\bf Theorem.}]\it
 }{\end{trivlist}\smallskip}
\newenvironment{theors}{\smallskip\begin{trivlist}
 \item[\hspace{\labelsep}{\noindent\bf Theorem}]\it
 }{\end{trivlist}\smallskip}
\newenvironment{lemma}{\smallskip\begin{trivlist}
 \item[\hspace{\labelsep}{\noindent\bf Lemma.}]\it
 }{\end{trivlist}\smallskip}
\newenvironment{lemmas}{\smallskip\begin{trivlist}
 \item[\hspace{\labelsep}{\noindent\bf Lemma}]\it
 }{\end{trivlist}\smallskip}
\newenvironment{propo}{\smallskip\begin{trivlist}
 \item[\hspace{\labelsep}{\noindent\bf Proposition.}]\it
 }{\end{trivlist}\smallskip}
\newenvironment{propos}{\smallskip\begin{trivlist}
 \item[\hspace{\labelsep}{\noindent\bf Proposition}]\it
 }{\end{trivlist}\smallskip}
\newenvironment{coro}{\smallskip\begin{trivlist}
 \item[\hspace{\labelsep}{\noindent\bf Corollary.}]\it
 }{\end{trivlist}\smallskip}
\newenvironment{defin}{\smallskip\begin{trivlist}
 \item[\hspace{\labelsep}{\noindent\em Definition.}]
 }{\end{trivlist}\smallskip}
\newenvironment{defins}{\smallskip\begin{trivlist}
 \item[\hspace{\labelsep}{\noindent\em Definition}]
 }{\end{trivlist}\smallskip}

\newenvironment{remark}{\begin{trivlist}
\item[\hspace{\labelsep}{\noindent\em Remark:}]
}{\end{trivlist}}
\newenvironment{proofs}{\begin{trivlist}  
\item[\hspace{\labelsep}{\noindent\em Proof:}]
}{\end{trivlist}}
\newcommand{\cuadro}{\hfill{$\qed$}}
\newcommand{\cuadrito}{\hbox{$\scriptstyle\sqcap \hskip-5.5pt \sqcup $}}

\def\lftcol#1{\vbox {\halign {##\hfil \cr #1\cr }}}
\long\def\direc#1#2{\hbox to \hsize{
     \lftcol{\hsize=7cm#1}\hfill\lftcol{\hsize=7cm#2}}}

\newcommand{\C}{\mathbb{C}}

\newcommand{\rD}{\mathrm{D}}
\newcommand{\h}{\mathbb{H}}
\newcommand{\M}{\mathbb{M}}
\newcommand{\N}{\mathbb{N}}
\newcommand{\p}{\mathbb{P}}

\newcommand{\bR}{\mathbf{R}}
\newcommand{\s}{\mathbb{S}}

\newcommand{\x}{\mathbb{X}}
\newcommand{\Z}{\mathbb{Z}}

\newcommand{\coh}{\hbox{\rm coh}}
\renewcommand{\deg}{\hbox{\rm deg}}
\newcommand{\degp}{\hbox{\tiny\rm deg}}

\renewcommand{\det}{\hbox{\rm det}}
\renewcommand{\dim}{\hbox{\rm dim}}

\newcommand{\Ext}{\hbox{\rm Ext}}
\newcommand{\End}{\hbox{\rm End}}

\newcommand{\Hom}{\hbox{\rm Hom}}
\newcommand{\id}{\hbox{\rm id}}

\renewcommand{\ker}{\hbox{\rm ker}}
\renewcommand{\max}{\hbox{\rm max}}
\renewcommand{\mod}{\hbox{\rm mod}}
\newcommand{\proj}{\mathrm{proj}}
\newcommand{\Mod}{\hbox{\rm Mod}}
\newcommand{\rad}{\hbox{\rm rad}}
\newcommand{\Root}{\hbox{\rm Root}}
\newcommand{\Spec}{\hbox{\rm Spec}}

\newcommand{\Oo}{\mathcal{O}}
\newcommand{\Tt}{\mathcal{T}}
\newcommand{\rperp}[1]{{#1^{\perp}}}
\newcommand{\lperp}[1]{{{}^{\perp}#1}}
\newcommand{\incl}{\hookrightarrow}
\newcommand{\ra}{\rightarrow}
\newcommand{\Mm}{\mathcal{M}}
\newcommand{\ZZ}{\mathbb{Z}}
\newcommand{\ZZplus}{{\mathbb{Z}_+}}
\newcommand{\la}{\lambda}
\newcommand{\si}{\sigma}
\newcommand{\vx}{\vec{x}}
\newcommand{\vc}{\vec{c}}
\newcommand{\vom}{\vec{\omega}}
\newcommand{\Lp}{{\mathbb{L}(p)}}
\newcommand{\HHplus}{\mathbb{H}_+}
\newcommand{\CC}{\mathbb{C}}
\newcommand{\Dsg}{\mathrm{D}_{\mathrm{Sg}}}
\newcommand{\Mmplus}{\mathcal{M}_+}
\newcommand{\Ppplus}{\mathcal{P}_+}
\newcommand{\Ssplus}{\mathcal{S}_+}
\newcommand{\Dd}{\mathcal{D}}
\newcommand{\Gaplus}{\Gamma_{+}}
\newcommand{\Dsing}{\mathrm{D}_{\mathrm{Sg}}}
\newcommand{\ga}{\gamma}
\newcommand{\al}{\alpha}
\newcommand{\be}{\beta}
\newcommand{\MCM}{\mathrm{MCM}}
\newcommand{\stabMCM}{\underline{\MCM}}
\newcommand{\hy}{\textrm{-}}
\newcommand{\bu}{\bullet}
\newcommand{\om}{\omega}

\newcommand{\euler}[1]{\chi_{#1}}
\def\raya{\raise1.5pt\hbox to 25pt{\vrule height1.5pt depth-1pt
           width25pt}}
\def\rayita{\raise2pt\hbox to 7.5pt{\vrule height1.5pt depth-1pt
           width7.5pt}}
\newcommand{\bulito}{\ {\scriptstyle \bullet}\, \ }
\newcommand{\bulitito}{{\scriptscriptstyle \bullet}}
\def\subsetnoteq{\mathbin{\hbox{$\subseteq \joinrel \hskip-8pt \lower3pt
                 \hbox{$\scriptscriptstyle /$}\ $}}}

\begin{document}

\title{Extended canonical algebras and Fuchsian singularities}
\author{Helmut Lenzing and Jos\'e A. de la Pe\~na}

\maketitle
\baselineskip=18pt plus 1 pt minus 1 pt

\section{Introduction} 

\subsection{Main results} 
\label{sect:1.1}
 Let $k$ be any field. Consider the hereditary
algebra $A=k\Delta$ associated to a finite connected quiver without
oriented cycles. A fundamental fact in Representation Theory of
Algebras is the distinction of the representation type of $A$: $A$
is {\em representation-finite\/} (that is, there are only finitely
many indecomposable $A$-modules, up to isomorphism) exactly when the
underlying graph $|\Delta|$ of $\Delta$ is of Dynkin type; $A$ is
{\em tame\/} (that is, for each $d\in \N$, the indecomposable
$d$-dimensional $A$-modules may be classified in a finite number of
one-parameter families of modules) exactly when $|\Delta|$ is of
extended Dynkin type; in the remaining cases $A$ is {\em wild\/}
(that is, there is a embedding $\mod\hy {k\langle x,y\rangle}\to
\mod\hy A$ which preserves indecomposability and isomorphism classes
from the category of finite dimensional modules over the ring in two
non-commuting indeterminates into the category of finite dimensional
$A$-modules). Consider the Auslander-Reiten translation $\tau_A$ in
$\mod\hy A$ and $P$ an indecomposable projective $A$-module. In case
$A$ is tame representation-infinite, the sequence of modules
$(\tau^{-n}_AP)_n$ is well defined and the algebra
$R(A,P)=\bigoplus\limits^\infty_{n=0}\Hom_A(P,\tau^{-n}_AP)$ is an
infinite dimensional positively $\Z$-graded surface singularity. The
algebra $R(A,P)$ reflects many properties of $A$ and $P$ and has a
particularly interesting structure, as shown in \cite{20}:

Assume $A=\C \tilde{\Delta}$ is a tame hereditary algebra where
$\tilde{\Delta}$ extends the Dynkin type $|\Delta|$. Let $P$ the
indecomposable projective associated to the vertex in
$\tilde{\Delta}\setminus \Delta$. Then $R(A,P)$ is isomorphic to the
algebra of invariants $\C [x,y]^G$, where $G\subset SL(2,\C)$ is a
binary polyhedral group of type $|\Delta|$. Accordingly the
completion of the graded algebra $R(A,P)$ is isomorphic to the
surface singularity of type $|\Delta|$.

In \cite{11}, Ringel introduced the canonical algebras
$C=C(p,\lambda)$ depending on a weight sequence $p=(p_1,\ldots,p_t)$
of positive integers and a parameter sequence $\lambda
=(\lambda_3,\ldots,\lambda_t)$ of pairwise distinct non-zero
elements from $k$. In \cite{5} it was shown that $\mod\hy C$ is
derived equivalent to $\coh\,(\x)$ the category of coherent sheaves
on a weighted projective line $\x=\x (p,\lambda)$.

One of the aims of this work is the introduction of a class of
algebras with related interesting properties. Let $P$ be an
indecomposable projective module over a canonical algebra
$C=C(p,\lambda)$, the one-point extension $A=C[P]$ defined as the
matrix algebra
\vglue-22truept
$$\left[\begin{matrix}k&0\\ P&C\end{matrix}\right]$$
is called an {\em extended canonical algebra}. In
section~\ref{sect:extended} we show that for two indecomposable
projective $C$-modules $P$ and $P'$, the algebras $C[P]$ and $C[P']$
are derived equivalent. Moreover, if $C$ is of tame type, the
extended canonical algebra is derived equivalent to a wild
hereditary or a wild canonical algebra, so essentially $C[P]$
belongs to a well-studied class of algebras. There are interesting
phenomena arising when $C$ (and hence $A$) is of wild type.

Consider the {\em Coxeter transformation\/} of $A$ as an automorphism
$\varphi_A\colon K_0(A)\to K_0(A)$ of the Grothendieck group of $A$,
given on the classes of indecomposable projective modules by the
formula $\varphi_A([P(S)])=-[I(S)]$, where $P(S)$ (resp. $I(S)$) is
the projective cover (resp. injective envelope) of a simple module
$S$. The characteristic polynomial $f_A(T)$ of $\varphi_A$ is called
the {\em Coxeter polynomial\/} of $A$.

Let $\euler{\x}=2-\sum_{i=1}^t(1-1/p_i)$ be the (orbifold) Euler
characteristic of $\x$. As shown in \cite{20}, for $\euler{\x}>0$
the classification problem of $\coh\,(\x)$ is related to the problem
of classifying the Cohen-Macaulay modules over a simple surface
singularity and is in fact equivalent to the problem of classifying
the graded Cohen-Macaulay modules over a corresponding
quasi-homogeneous singularity. Assume $\euler{\x}<0$. By \cite{2,22}
we know that for $A=C[P]$ an extended canonical algebra
$$f_A(T)=P_C(T)f_C(T)$$
where $P_C(T)$ is the Hilbert-Poincar\'e series of the positively
graded algebra
$R(p,\lambda)=\bigoplus\limits^\infty_{n=0}\Hom_C(M,\tau^n_CM)$,
where $M$ is a rank one not preprojective $C$-module. Equivalently,
$R(p,\lambda)=\bigoplus\limits^\infty_{n=0}\Hom(\Oo,\tau_\x^n\Oo)$,
where $\Oo$ is the structure sheaf on $\x$. Recall from \cite{1,2}
that in case $k=\C$, we can interpret $R(p,\lambda)$ as an algebra
of entire automorphic forms associated to the action of a suitable
{\em Fuchsian group\/} of the first kind, acting on the upper half
plane $\h_+$.

>From \cite{1}, we know that the $k$-algebra $R=R(p,\lambda)$ is
commutative, graded integral Gorenstein, in particular
Cohen-Macaulay, of Krull dimension two. The complexity of the
surface singularity $R$ is described by the triangulated category
$$
\Dsing^\ZZ(R)=\frac{\rD^b(\mod^\ZZ\hy R)}{\rD^b(\proj^\ZZ\hy R)},
$$
where $\mod^\ZZ\hy R$ (resp.\ $\proj^\ZZ\hy R$)  denotes the
category of finitely generated (resp.\ finitely generated
projective) $\ZZ$-graded $R$-modules. This category was considered
by Buchweitz \cite{36} and Orlov \cite{34}, see also Krause's
account \cite{37} for a related, but slightly different approach.
For $\euler{\x}>0$, where the weight type of $\x$ determines a
Dynkin quiver $\Delta$, Kajiura, Saito, Takahashi and Ueda \cite{35}
have shown that $\Dsing^\ZZ(R)$ is equivalent to the derived
category of finite dimensional modules over the path algebra
$k\,\Delta$. For $\euler{\x}=0$, the algebra $R$ hence
$\Dsing^\ZZ(R)$ is not defined, but a close variant
$\Dsing^{\Lp}(S)$, as shown by Ueda~\cite{38}, is equivalent to the
derived category $\rD^b(\coh\,(\x))$ of coherent sheaves on a
weighted projective line which is tubular, that is, has weight type
$(2,3,6)$, $(2,4,4)$, $(3,3,3)$ or $(2,2,2,2)$.

In section~\ref{extended:canonical} we deal with the case
$\euler{\x}<0$ and prove that this category, as first observed by
Saito and Takahashi (for the field $\C$ of complex numbers), is
described in the following way.

\begin{theorem} \label{thm:1}
Let $k$ be an algebraically closed field. Assume $\euler{\x}<0$ and
let $R$ be the positively $\ZZ$-graded surface singularity attached
to $\x$. Then there exists a tilting object $\bar{T}$ in the
triangulated category $\Dsing^\ZZ(R)$ whose endomorphism ring is
isomorphic to an extended canonical algebra $C[P]$, where $C$ is the
canonical algebra associated with $\x$.
\end{theorem}
It follows that the categories $\Dsing^\ZZ(R)=\Tt$ and
$\rD^b(\mod(C[P]))$ are equivalent as triangulated categories. In
section~\ref{coxeter:dynkin} we further introduce the concept of
Coxeter-Dynkin algebras and establish their relationship to the
Coxeter-Dynkin diagrams from singularity theory.

More precise information on the structure of the ring $R(p,\la)$ is
obtained by a closer examination of the spectral properties of the
Coxeter transformation of the extended canonical algebra $C[P]$.

A sequence of weights $p=(p_1,\ldots,p_t)$ will always satisfy
$p_1\le p_2\le \cdots \le p_t$. We consider the lexicographical
ordering of sequences $(p_1,\ldots,p_t)\le (q_1,\ldots,q_s)$ if
$t=s$ and $p_i\le q_i$ for $1\le i\le t$. Extend the relation $p\le
q$ (and say that $q$ {\em dominates\/} $p$) to weight sequences of
(possibly) different length by adding $1$'s if necessary. The
following result is shown in section~\ref{sect:cox} based on
techniques developed in \cite{3} and will be fundamental in the
proof of the main results.

\begin{theorem}\label{thm:2}
Let $A=C[P]$ be an extended canonical algebra of the \emph{wild}
canonical algebra $C=C(p,\lambda)$. The following happens:

\begin{itemize}
\item[{\rm (a)}] $f_A(T)$ has at most $4$ roots not in $\s^1$.

\item[{\rm (b)}] The roots of $f_A(T)$ lie on the unit circle $\s^1$
if and only if the weight sequence $p$ belongs to the 38-member list
determined by all $p<q$ with $q$ belonging to the following {\em
critical list\/}:

\begin{itemize}
\item[{\rm ($t=3$):}] $(2,3,11)$, $(2,4,9)$, $(2,5,8)$, $(2,6,7)$,\\
$(3,3,8)$, $(3,4,7)$, $(3,5,6)$,\\ $(4,4,6)$, $(4,5,5)$,

\item[{\rm ($t=4$):}] $(2,2,2,7)$, $(2,2,3,6)$, $(2,3,4,4)$,
$(3,3,3,4)$,

\item[{\rm ($t=5$):}] $(2,2,2,2,5)$, $(2,2,2,3,4)$, $(2,2,3,3,3)$,

\item[{\rm ($t=6$):}] $(2,2,2,2,2,3)$,

\item[{\rm ($t=7$):}] $(2,2,2,2,2,2,2)$.
\end{itemize}
\end{itemize}
\end{theorem}

We shall say that the algebra $R(p,\lambda)$ (and also the weight
sequence $p$) is {\em formally $n$-generated\/} if
$$P_C(T)=\frac{\prod\limits^{n-2}_{i=1}
(1-T^{c_i})}{\prod\limits^n_{j=1}(1-T^{d_j})}$$ for certain natural
numbers $c_1,\ldots,c_{n-2}$ and $d_1,\ldots,d_n$, all $\ge 2$. The
algebra $R(p,\lambda)$ (and also the weight sequence $p$) is {\em
formally a complete intersection\/} if $P_C(T)$ is a rational
function $f_1(T)/f_2(T)$, where each $f_i(T)$ is a product of
cyclotomic polynomials.

\begin{theorem}\label{thm:3}
Let $C=C(p,\lambda)$ be a wild canonical algebra with weight
sequence $(p_1,\ldots,p_t)$ and $A=C[P]$ be an extended canonical
algebra. The following are equivalent:

\begin{itemize}
\item[{\rm (a)}] $R(p,\lambda)$ is formally $3$- or $4$-generated

\item[{\rm (b)}] $R(p,\lambda)$ is formally a complete intersection

\item[{\rm (c)}] The roots of $f_A(T)$ lie on $\s^1$.
\end{itemize}

Moreover, for $t=3$ the algebra $R(p,\lambda)$ is a graded complete
intersection of the form $k[X_1,\ldots,X_s]/(\rho_3,\ldots,\rho_s)$
where $s=3$ or $4$ and $\rho_3,\ldots,\rho_s$ is a homogeneous
regular sequence. For $k=\C$ the assertion also holds for $t\geq4$
for $R(p,\la')$ for a suitable choice of parameters.
$\la'=(\la'_3,\ldots,\la'_t)$.
\end{theorem}

We remark that in almost every case $\Root\,f_A(T)\subset \s^1$
implies that the Coxeter transformation $\varphi_A$ is periodic. In
fact the weight sequences $(3,3,3,3)$ and $(2,2,2,2,4)$ are the only
exceptions (section~\ref{sect:cox}).

\begin{theorem}\label{thm:4}
Let $C=C(p,\lambda)$ be a wild canonical algebra with weight
sequence $p=(p_1,\ldots,p_t)$. Consider $A=C[P]$ an extended
canonical algebra. The following are equivalent:

\begin{itemize}
\item[{\rm (a)}] $R(p,\lambda)$ is formally $3$-generated

\item[{\rm (b)}] $\varphi_A$ is periodic of period $d$ and there is a
primitive $d$-th root of unity which is root of $f_A(T)$.
\end{itemize}

For $t=3$ the algebra $R(p,\lambda)$ is always a graded complete
intersection of the form $k[X_1,X_2,X_3]/(f)$. Moreover, for
$t\geq4$ and $k=\C$ this also holds for $R(p,\lambda')$ for a
suitable choice of parameters $\la'=(\la'_3,\ldots,\la'_t)$.
\end{theorem}

For the proof of Theorem~\ref{thm:3} (resp. Theorem~\ref{thm:4}) we
classify in section~\ref{sect:poinc} all the weight sequences $p$
such that $R(p,\lambda)$ is formally a complete intersection (resp.
$R(p,\lambda)$ has $3$ homogeneous generators). In the complex case
the algebras $R(p,\lambda)$ of Theorem~\ref{thm:3} correspond to the
Fuchsian singularities which are minimal elliptic
\cite[Proposition~5.5.1]{14} and the classification is related to
Laufer's \cite{15}. The algebras $R(p,\lambda)$ of
Theorem~\ref{thm:4} relate to classifications by Dolgachev \cite{7}
and Wagreich \cite{8} and include the $14$ exceptional unimodal
Arnold's singularities \cite{25}. We refer the reader to the
complete account by Ebeling \cite{10}.

The research for this work was done during exchange visits
M\'exico-Paderborn. We thank our universities and CONACyT in
M\'exico for support. We thank Henning Krause for directing our
attention to Buchweitz's paper~\cite{36}.

\subsection{Notation and conventions} 
Let $Q$ be a finite quiver without oriented cycles. The {\em path
algebra\/} $kQ$ has as basis all the oriented paths in $Q$ and
product given by juxtaposition of paths. Given an ideal $I$ of $kQ$
which is {\em admissible\/} (that is, $(kQ^+)^m\subset I\subset
(kQ^+)^2$ for some $m\ge 2$, where $kQ^+$ is the ideal of $kQ$
generated by the arrows), we consider the finite dimensional
$k$-algebra $A=kQ/I$. By `module' we mean a finite dimensional right
$A$-module. The category of modules is denoted $\mod\hy A$. A module
is identified with a covariant functor $X\colon kQ\to \mod\hy k$
such that $X(\rho)=0$ for every $\rho \in I$. Important modules are
the simple modules $S_i$ associated to vertices $i\in Q_0$ (= set of
vertices of $Q$), the projective cover $P_i=A(-,i)$ of $S_i$ and the
injective envelope $I_i=DA(i,-)$ of $S_i$, where $D=\Hom_k(-,k)$ is
the canonical duality.

We denote by $K_0(A)$ the Grothendieck group of $A$. Since $A$ has
finite global dimension, the classes $[P_i]$ with $i\in Q_0$ form a
basis of $K_0(A)$. Thus the {\em Coxeter transformation\/}
$\varphi_A\colon K_0(A)\to K_0(A)$, given by
$\varphi_A([P_i])=-[I_i]$ defines an isomorphism. The Grothendieck
group $K_0(A)$ is equipped with a bilinear form $\langle
-,-\rangle_A\colon K_0(A)\times K_0(A)\to \Z$, called the {\em Euler
form}, defined in the classes of modules $X$ and $Y$ as $\langle
[X],[Y]\rangle_A=\sum\limits^\infty_{i=0}(-1)^i\dim_k\Ext^i_A(X,Y)$.
In case $A=k\Delta$ is a hereditary algebra, then
$\varphi_A([X])=[\tau_AX]$ for any indecomposable non-projective
$A$-module $X$. For the general situation, we have to look at the
derived category $D(A)=\rD^b(\mod\hy A)$ of bounded complexes of
$A$-modules.

The {\em derived category\/} $D(A)$ contains a copy $\mod\hy A[n]$
of $\mod\hy A$ for each integer $n\in \Z$, with objects written
$X[n]$ and satisfying
$$\Hom_{D(A)}(X[n],Y[m])=\Ext^{n-m}_A(X,Y).$$
We say that an algebra $A$ is \emph{derived hereditary} (resp.\
\emph{derived canonical}) if $D(A)$ is triangle equivalent to $D(H)$
(resp.\ $D(C)$) for a hereditary algebra $H$ (resp.\ a canonical
algebra $C$).

The category $D(A)$ has Auslander-Reiten triangles which yield a
self-equivalence $\tau_{D(A)}$ of $D(A)$, the Auslander-Reiten
translation, satisfying
$\Hom_{D(A)}(Y,\tau_{D(A)}X[1])=D\Hom_A(X,Y)$. The natural
isomorphism $K_0(A)\to K_0(D(A))$, $X\mapsto X[0]$ yields
$\varphi_A([X])=[\tau_{D(A)}X]$.

For background material on representations of algebras and derived
categories we refer the reader to \cite{16,11}.

For vectors $v,w\in K_0(A)$ we get $\langle
v,\varphi_A(w)\rangle_A=-\langle w,v\rangle_A$.

\section{Extended canonical algebras: basic properties} 
\label{sect:extended}
\subsection{} 
Let $C=C(p,\lambda)$ be the {\em canonical algebra\/} defined by the
weight sequence $p=(p_1,\ldots,p_t)$ with $p_i\ge 2$ and $(\lambda
=\lambda_3,\ldots,\lambda_t)$ a sequence of pairwise distinct
non-zero elements of $k$, that is, $C$ is defined by the quiver

{\small
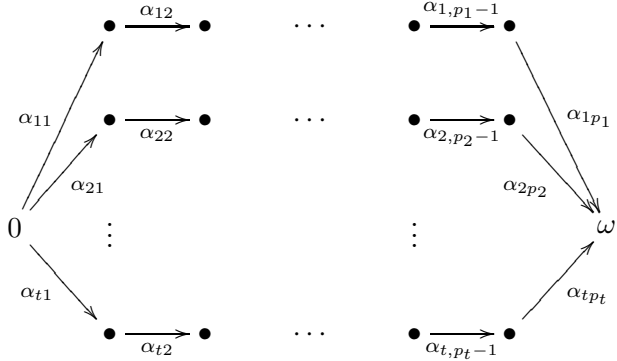
$$
\xymatrix{
                                                                  &\bu\ar[r]^{\al_{12}}&\bu    &\cdots&\bu\ar[r]^{\al_{1,p_1-1}}&\bu\ar[ddr]^{\al_{1p_1}}& &\\
                                                                  &\bu\ar[r]_{\al_{22}}&\bu    &\cdots&\bu\ar[r]_{\al_{2,p_2-1}}&\bu\ar[dr]_{\al_{2p_2}}&  &\\
 0\ar[ur]_{\al_{21}}\ar[uur]^{\al_{11}}\ar[dr]_{\al_{t1}}         &\vdots              &       &      &\vdots                   &                       &\om&\\
                                                                  &\bu\ar[r]_{\al_{t2}}&\bu
                                                                  &\cdots&\bu\ar[r]_{\al_{t,p_t-1}}&\bu\ar[ur]_{\al_{tp_t}}&
                                                                  &
}
$$}
\vglue-12truept \noindent satisfying the $t-2$ equations
$$\alpha_{ip_i}\ldots \alpha_{i2}\alpha_{i1}=\alpha_{2p_2}\ldots
\alpha_{22}\alpha_{21}-\lambda_i\alpha_{1p_1}\ldots
\alpha_{12}\alpha_{11},\ i=3,\ldots,t.$$

The algebra $C$ is a one-point extension $H[M]$ of the hereditary
algebra $H=C/(0)$ by an $H$-module $M$ with dimension vector
$[M]=(\dim_kM(i))_i\in K_0(H)$ as follows:

{\tiny
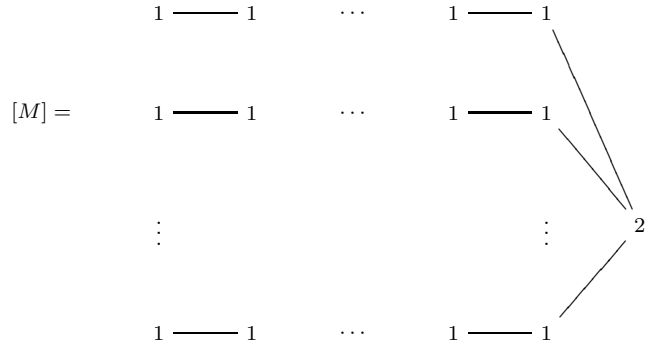
$$
\xymatrix{
    &1\ar@{-}[r]&1&\cdots&1\ar@{-}[r]&1\ar@{-}[ddr]&    \\
[M]= &1\ar@{-}[r]&1&\cdots&1\ar@{-}[r]&1\ar@{-}[dr] &    \\
     &\vdots     & &      &&\vdots    &  2 \\
    &1\ar@{-}[r]&1&\cdots&1\ar@{-}[r]&1\ar@{-}[ur]&
 }
$$}
 \noindent In case $t\ge 3$, the module $M$ is an
indecomposable $H$-module which is not preprojective or
preinjective.

Observe that the underlying graph of $H$ is a star
$[p_1,p_2,\ldots,p_t]$ with linear arms having $p_i$ vertices,
$i=1,\ldots,t$. If $H$ is representation-finite, then
$[p_1,p_2,\ldots,p_t]$ is a Dynkin diagram (that is,
$\sum\limits^t_{i=1}\frac{1}{p_i}>t-2$) and $C$ is {\em tame of
domestic type}. If $H$ is tame, then $[p_1,p_2,\ldots,p_t]$ is an
extended Dynkin diagram (that is,
$\sum\limits^t_{i=1}\frac{1}{p_i}=t-2$) and $C$ is {\em tame of
tubular type}. See \cite{11} for details.

\subsection{} 
The representation theory of $\mod\hy C$ for $C=C(p,\lambda)$ a
canonical algebra is controlled by the category $\coh\,(\x)$ of
coherent sheaves on a weighted projective line $\x =\x (p,\lambda)$,
since the derived categories $D(C)$ and $\rD^b(\coh\,(\x))$ are
equivalent as triangulated categories \cite{20}. The complexity of
the classification problem for $\coh\,(\x)$, and hence for $\mod\hy
C$, is essentially determined by the (orbifold) Euler characteristic
$\euler{\x}=2-\sum_{i=0}^t(1-1/p_i)$. Indeed, for $\euler{\x}>0$,
the algebra $C$ is tame of domestic type and for $\euler{\x}=0$, the
algebra $C$ is tubular. The wild case $\euler{\x}<0$ was carefully
studied in \cite{2}, a paper which is the basis of the present
investigation.

\subsection{} 
Let $C=C(p,\lambda)$ be a canonical algebra. Let $P$ be an
indecomposable projective or injective $C$-module, then $A=C[P]$ is
called an {\em extended canonical\/} algebra. Hence $A$ arises from
by adjoining one arrow with a new vertex to an arbitrary vertex of
$C$ while keeping the relations for $C$ without introducing any new
relations. In particular, the opposite algebra of an extended
canonical algebra is again extended canonical.

\begin{lemma}
Any extended canonical algebra is wild.
\end{lemma}

\begin{proofs}
Recall that for an algebra $B=kQ/I$ where $Q$ has no oriented cycles
and $I$ is generated by $\rho_1,\ldots,\rho_s\in
\bigcup\limits_{i,j\in Q_0}I(i,j)$, the {\em Tits\/} (quadratic) {\em
form\/} $q_B\colon K_0(B)\to \Z$ is defined by
$${\textstyle q_B(x)=\sum\limits_{i\in Q_0}x(i)^2-\sum\limits_{i\to
j}x(i)x(j)+\sum\limits_{i,j\in Q_0}r(i,j)x(i)x(j),}$$
where $r(i,j)=\# \{s\colon \rho_s\in I(i,j)\}$. The Tits form is
weakly non-negative (i.e. $q_A(v)\ge 0$ for $v\in \N^{Q_0}$) if $B$
is a tame algebra \cite{17}.

For an extended canonical algebra $A=C[P]$ with extension vertex $*$
such that $\rad\,P_*=P=P_j$ for some vertex $j$ in $C$, we have the
following:

\begin{itemize}
\item[{\rm $\bulito$}] since $gl\,\dim\,A=2$, then $q_A(x)=\langle
x,x\rangle_A$;

\item[{\rm $\bulito$}] the vector $w\in K_0(C)\subset K_0(A)$ with
$w(i)=1$ for every $i$ in $C$, satisfies $q_c(w)=0$;

\item[{\rm $\bulito$}] for $e_*=[S_*]$, we get
$$q_A(2w+e_*)=4q_c(w)-2w(j)+1<0.$$
Hence $A$ is of wild type.\cuadro
\end{itemize}
\end{proofs}

\subsection{} 
\label{sect:2.4} The following result is fundamental for introducing
the concept of extended canonical algebras.

\begin{propo}
Let $X$ and $Y$ be two indecomposable projective or injective
$C$-modules over a canonical algebra $C$. Then the extended
canonical algebras $C[X]$ and $C[Y]$ are derived equivalent.
\end{propo}
In particular, the derived class of an extended canonical algebra is
independent of the chosen projective module.

\begin{proof}
(See \cite{3}). Under the equivalence $D(C)=\rD^b(\coh\,(\x
(p,\lambda)))$ the modules $X$ and $Y$ become line bundles over $\x$
(up to translation in $\rD^b(\coh\,(\x))$). Since the Picard group
of $\x$ acts transitively on isomorphism classes of line bundles,
there is a self-equivalence of $D(C)$ sending $X$ to $Y$. The
assertion follows from \cite{13}.
\end{proof}

\subsection{} 
The following remark is useful: \label{sect:2.5}
\begin{propo}
Let $A$ be an extended canonical algebra of $C=C(p,\lambda)$ with
$p=(p_1,\ldots,p_t)$ satisfying $t\ge 3$. Then $A$ is derived
equivalent to a one-point extension $H[N]$ of a hereditary algebra
$H$ by an indecomposable module $N$.
\end{propo}

\begin{proof}
By (2.3), we may assume that $A$ is the path algebra of the quiver

{\small
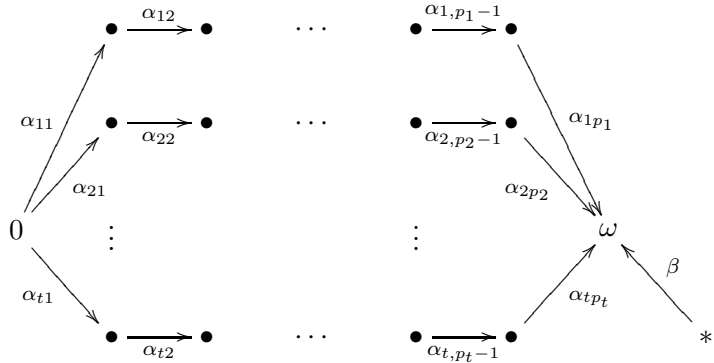
$$
\xymatrix{
                                                                  &\bu\ar[r]^{\al_{12}}&\bu    &\cdots&\bu\ar[r]^{\al_{1,p_1-1}}&\bu\ar[ddr]^{\al_{1p_1}}& \\
                                                                  &\bu\ar[r]_{\al_{22}}&\bu    &\cdots&\bu\ar[r]_{\al_{2,p_2-1}}&\bu\ar[dr]_{\al_{2p_2}}& \\
 0\ar[ur]_{\al_{21}}\ar[uur]^{\al_{11}}\ar[dr]_{\al_{t1}}         &\vdots              &       &      &\vdots                   &                       &\om\\
                                                                  &\bu\ar[r]_{\al_{t2}}&\bu
                                                                  &\cdots&\bu\ar[r]_{\al_{t,p_t-1}}&\bu\ar[ur]_{\al_{tp_t}}&&\ast\ar[lu]_{\be}\\
}
$$}

\noindent equipped with the canonical relations $\alpha_{ip_i}\ldots
\alpha_{i1}=\alpha_{2p_2}\ldots
\alpha_{21}-\lambda_i\alpha_{1p_1}\ldots \alpha_{11}$ ($3\le i\le
t$). It follows that $A$ is the one-point extension of the path
algebra of the star\break $[2,p_1,p_2,\ldots,p_t]$ by an
indecomposable module $N=N(\lambda_3,\ldots,\lambda_t)$ whose
restriction to $[p_1,p_2,\ldots,p_t]$ is $M$ as in (2.1) and $N(*)=0$.
\end{proof}

\begin{remark}
An obvious variant of the above statement yields families of
one-point extensions of hereditary algebras which are pairwise
derived equivalent. For instance, for the canonical type $(2,3,7)$,
the Proposition yields $10=(2-1)+(3-1)+(7-1)+1$ choices of pairs
$(H,N)$ of a hereditary algebra $H$ and an indecomposable $H$-module
$N$ such that $H[N]$ is derived equivalent to an extended canonical
algebra of type $(2,3,7)$.
\end{remark}

\section{The derived category of an extended canonical algebra} 
\label{extended:canonical}
\subsection{} In this section we are going to investigate the nature of
the bounded derived category of an extended canonical algebra
$A=C[P]$. As it turns out the structure of this triangulated
category will sensibly depend on the sign of the (orbifold) Euler
characteristic $\euler{\x}$ of the weighted projective line $\x$
associated to $C$.

Let $\Tt$ be a triangulated $k$-category, see \cite{26,27,28} for
definition and properties. An object $E$ in $\Tt$ is called
\emph{exceptional} if $\End(E)=k$ and $\Hom(E,E[n])=0$ for all
integers $n\neq0$. By $\lperp{E}$ (resp.\ $\rperp{E}$ we denote the
full triangulated subcategory of $\Tt$ consisting of all objects
$X\in\Tt$ (resp.\ $Y\in\Tt$) satisfying $\Hom(X,E[n])=0$ (resp.\
$\Hom(E[n],Y)=0$ for each integer $n$. By \cite{29} the inclusion
$\lperp{E}\incl\Tt$ (resp.\ $\rperp{E}\incl\Tt$) admits an exact
left (resp.\ right) adjoint $\ell:\Tt\ra \lperp{E}$ (resp.\
$r:\Tt\ra\rperp{E}$).

For the purpose of this paper we call an exceptional object
\emph{special in $\Tt$} if one of the following two conditions is
satisfied:

$(i)$ the left perpendicular category $\lperp{E}$ is equivalent to
$\rD^b(\coh\, (\x))$ for some weighted projective line $\x$ and,
moreover, the left adjoint $\ell$ maps $E$ to a line bundle in
$\coh(\x)$.

$(ii)$ the right perpendicular category $\rperp{E}$ is equivalent to
$\rD^b(\coh\,(\x))$ for some weighted projective line $\x$ and,
moreover, the right adjoint $r$ maps $E$ to a line bundle in
$\coh(\x)$.

Again, for the purpose of this paper, an object $T$ of $\Tt$ is
called a \emph{tilting object} in $\Tt$ if $(i)$ $T$ generates $\Tt$
as a triangulated category, $(ii)$ $\Hom(T,T[n])=0$ holds for each
non-zero integer $n$, and $(iii)$ the endomorphism algebra of $T$
has finite global dimension.

\subsection{} Our main tool to investigate the shape of $D(A)$ is the
following proposition. \label{propo:special}
\begin{propo}
Let $\Tt$ be a triangulated category having an exceptional object
$E$ that is special in $\Tt$. Then there exists a tilting object
$\bar{T}$ of $\Tt$ whose endomorphism ring is an extended canonical
algebra.
\end{propo}

\begin{proof}
With the previous notations we assume that
$\rperp{E}=\rD^b(\coh(\x))$ and $r(E)$ is a line bundle in
$\coh\,(\x)$. (The assumption on the left perpendicular $\lperp{E}$
category is treated similarly.)  We choose a tilting object $T$ in
$\coh\,(\x)$, hence in $\rD^b(\coh\,(\x))$ having the line bundle
$r(E)$ as a direct summand and such that $\End(T)=C$ is the
canonical algebra attached to $\x$ (see \cite{5}). We claim that
$\bar{T}=T\oplus E$ is a tilting object in $\Tt$. Indeed, since $T$
generates $\rD^b(\coh\,\x)$ and $E$ together with $\rperp{E}$
generates $\Tt$, it follows that $\bar{T}$ generates $\Tt$. Next, we
show that $\Hom(\bar{T},\bar{T}[n])=0$ for each nonzero integer $n$.
This reduces to show that $\Hom(E[n],T)=0$ and $\Hom(T,E[n])$ holds
for every nonzero $n$. The first assertion holds for each $n$ since
$T$ belongs to $\rperp{E}$. For the second we use that
$\Hom(T,E[n])=\Hom(T[-n],r(E))$ is zero since by construction $r(E)$
is a direct summand of the tilting object $T$. Finally, the
endomorphism ring of $\bar{T}$ is given as the matrix ring
$$
\left(
\begin{array}{cc}
 C & 0 \\
 P  & k
 \end{array}
 \right), \quad \textrm{ where } P=\Hom(T,E)=\Hom(T,rE)
$$
is an indecomposable projective $C$-module, hence
$\End(\bar{T})=C[P]$ is an extended canonical algebra. Moreover, as
it is easily seen, $C[P]$ has global dimension two.
\end{proof}

Keeping the assumptions on $E$ and $\Tt$ from the proposition, we
obtain. \label{coro:algebraic}
\begin{coro}
If $\Tt$ is triangle equivalent to a bounded derived category $D(B)$
for some finite dimensional $k$-algebra $B$ or, more generally, if
$\Tt$ is algebraic in the sense of Keller \cite{28}, then $\Tt$ is
triangle equivalent to $\rD^b(\mod\hy A)$, for the extended
canonical algebra $A=C[P]$.
\end{coro}
\begin{proof}
The first claim follows from \cite{30}, the general version requires
\cite{28} or \cite{31}.
\end{proof}

\subsection{Positive Euler characteristic: the domestic case} 
Consider a canonical algebra $C=C(p_1,\ldots,p_t)$ of domestic type,
that is, $\sum\limits^t_{i=1}\frac{1}{p_i}>1$. Let $\Delta
=[p_1,\ldots,p_t]$ be the star corresponding to the weight sequence
$p=(p_1,\ldots,p_t)$ and $\tilde{\Delta}$ be the corresponding
extended Dynkin diagram. Then $\tilde{\Delta}$ admits a unique
positive {\em additive function\/} $\lambda$ assuming value $1$, that
is, $\lambda \colon \tilde{\Delta}\to \N$ satisfies the conditions:

(i) $2\lambda (i)=\sum\limits_{j\in i^+}\lambda (j)$, where $i^+$ is
the set of neighbors of $i$;

(ii) for some vertex $i\in \tilde{\Delta}_0$, $\lambda (i)=1$. Such
an $i$ is called an {\em extension vertex}.

The {\em double extended graph\/} of type $\Delta$, denoted by
$\tilde{\tilde{\Delta}}$, is the graph arising form $\tilde{\Delta}$
by adjoining a new edge in an extension vertex. The list of double
extended Dynkin graphs is the following:

 \def\entrya{$\scriptscriptstyle
\xymatrix{
                                    &\bu\ar@{-}[r] &\bu&\cdots&\bu\ar@{-}[r]&\bu\ar@{-}[dr]  &            &   \\
 \bu\ar@{-}[ur]\ar@{-}[dr]&         &              &          &            &             &\bu\ar@{-}[r]&\ast \\
                                    &\bu\ar@{-}[r] &\bu&\cdots&\bu\ar@{-}[r]&\bu\ar@{-}[ur]  &            &
} $}
 \def\entryb{$\scriptscriptstyle
 \xymatrix{
 \bu\ar@{-}[dr] &             &    &       &             &              &\bu                    &   \\
                &\bu\ar@{-}[r]&\bu &\cdots                &\bu\ar@{-}[r]&\bu\ar@{-}[ur]\ar@{-}[dr]&   \\
 \bu\ar@{-}[ur] &             &    &       &             &              &\bu\ar@{-}[r]           &\ast
 }
$}
 \def\entryc{$\scriptscriptstyle
\xymatrix{
 \bu\ar@{-}[r]& \bu\ar@{-}[r]& \bu\ar@{-}[r]\ar@{-}[d] & \bu\ar@{-}[r] & \bu\ar@{-}[r] & \ast & (=[3,3,4]) \\
              &              & \bu\ar@{-}[d]           &               &               &      &     \\
              &              & \bu                     &               &               &      &
 }$}
 \def\entryd{$\scriptscriptstyle
 \xymatrix{
 \bu\ar@{-}[r]& \bu\ar@{-}[r]& \bu\ar@{-}[r]& \bu\ar@{-}[r]\ar@{-}[d]& \bu\ar@{-}[r]& \bu\ar@{-}[r]& \bu\ar@{-}[r]& \ast&(=[2,4,5])\\
              &              &              & \bu                    &              &              &              & &} $}
 \def\entryf{$\scriptscriptstyle
 \xymatrix{
 \bu\ar@{-}[r] & \bu\ar@{-}[r] &\bu\ar@{-}[r]\ar@{-}[d] & \bu\ar@{-}[r]& \bu\ar@{-}[r]&\bu\ar@{-}[r]& \bu\ar@{-}[r] &\bu\ar@{-}[r] &\ast &(=[2,3,7]) \\
               &               & \bu                    &              &              &             &               &              &     & }$}

{\scriptsize
\begin{center}
\begin{tabular}{l|l}
weight\\ sequence &\ \ double extended Dynkin graph\\
\hline
&\\
$(p,q)$ &\ \ \entrya\\ 
 & \\
$(2,2,4)$ &\ \ \entryb\\ 
 & \\
$(2,3,3)$ &\ \ \entryc \\ 
 & \\
$(2,3,4)$ &\ \ \entryd\\ 
 & \\
$(2,3,5)$ &\ \ \entryf
\end{tabular}
\end{center}}
\label{propo:hereditary}
\begin{propo}
Let $C=C(p,\lambda)$ be a canonical algebra of domestic type
$p=(p_1,p_2,p_3)$. Let $\Delta =[p_1,p_2,p_3]$ be the associated
Dynkin diagram. For any indecomposable preprojective $C$-module $N$,
the one-point extension $A=C[N]$ is derived equivalent to a
hereditary algebra of type $\tilde{\tilde{\Delta}}$. In particular,
an extended canonical algebra of weight type $p$ is derived
hereditary of type $\tilde{\tilde{\Delta}}$.
\end{propo}

\begin{proof}
The algebra $C$ is tilted of a hereditary algebra
$H=k\tilde{\Delta}_1$, where $\tilde{\Delta}_1$ is a quiver with
underlying graph $\tilde{\Delta}$. There is a derived equivalence
$F\colon \rD^b(\mod\hy H)\to \rD^b(\mod\hy C)$ sending the
indecomposable projective $P_i$ corresponding to an extension vertex
$i$ of $\tilde{\Delta}$ into the projective $C$-module $P_\omega$.
Observe that $H[P_i]=k\tilde{\tilde{\Delta}}_1$ is a hereditary
algebra where $\tilde{\tilde{\Delta}}_1$ is a quiver with underlying
graph $\tilde{\tilde{\Delta}}$.

Let $\x =\x (p,\lambda)$ be a weighted projective line such that
$\rD^b(\mod\hy C)=\rD^b(\coh\,(\x))$. As an object in $\coh\,(\x)$
the object $P_\omega$ has rank one (see \cite{5}). Also by \cite{5},
there is an equivalence in $\rD^b(\coh\,(\x))$ sending $P_\omega$ to
any indecomposable preprojective $C$-module $N$. By \cite{13}, the
one-point extension $C[N]$ is derived equivalent to $C[P_\omega]$
which is derived hereditary of type $\tilde{\tilde{\Delta}}$.
\end{proof}

\subsection{} 
The converse of Proposition~\ref{propo:hereditary} also holds.
\begin{propo}
Let $C=C(p,\lambda)$ be a canonical algebra and $P$ be an
indecomposable projective $C$-module. The extended canonical algebra
$A=C[P]$ is derived hereditary if and only if $C$ is tame domestic.
\end{propo}

\begin{proof}
If $C$ is tame domestic, then $A$ is derived hereditary by
(\ref{sect:2.5}). For the converse, consider the set of weight
sequences $p=(p_1,p_2,\ldots,p_t)$ with $2\le p_1\le p_2\le \cdots
\le p_t$ with the domination order defined in (\ref{sect:1.1}).

The statement follows by induction on the domination order from the
following two facts:

\begin{itemize}
\item[{\rm (a)}] a canonical tubular algebra $C$ is not derived
hereditary;

\item[{\rm (b)}] any wild weight sequence dominates a tubular one;

\item[{\rm (c)}] if $M$ is an indecomposable $B$-module such that
$B[M]$ is derived hereditary, then $B$ is derived hereditary.
\end{itemize}

(a): follows from the structure of derived categories of hereditary
algebras, see \cite{16}. (b) is clear.

(c): Assume $\rD^b(\mod\hy {B[M]})=\rD^b(\mod\,H)$ for a hereditary
algebra $H$. By \cite{20}, $\rD^b(\mod\hy B)$ is equivalent to the
right perpendicular category in $\rD^b(\mod\hy H)$ with respect to
an exceptional object $E$, that is, $E$ is an indecomposable object
satisfying $\Ext^1(E,E)=0$ and
$$\rD^b(\mod\hy B)=E^\perp =\{X\in \rD^b(\mod\hy H)\colon
\Hom\,(E,X)=0=\Ext^1(E,X)\}.$$ Without loss of generality we may
assume that $E\in \mod\hy H$. Then $\rD^b(\mod\hy B)\cong
\rD^b(E^\perp)$, where now $E^\perp$ is formed in $\mod\hy H$. Hence
$E^\perp =\mod\hy {H'}$ for a hereditary algebra $H'$.
\end{proof}

\subsection{Euler characteristic zero: the tubular case} 
Consider a canonical algebra $C=C(p,\lambda)$ with weight sequence
$p=(p_1,\ldots,p_t)$, we shall assume that $2\le p_1\le p_2\le
\cdots \le p_t$. The module category $\mod\hy C$ accepts a {\em
separating tubular family\/} $\mathcal{T}=(T_\lambda)_{\lambda \in
\p_1k}$, where $T_\lambda$ is a homogeneous tube for all $\lambda$
with the exception of $t$ tubes $T_{\lambda_1},\ldots,T_{\lambda_t}$
with $T_{\lambda_i}$ of rank $p_i$ ($1\le i\le t$). See \cite{11}.

Let $\x=\x(p,\lambda)$ be the weighted projective line such that
$\mod\hy C$ and $\coh\,(\x)$ are derived equivalent. We fix an
equivalence $\rD^b(\mod\hy C)=\rD^b(\coh\,(\x))$. Let $S$ be a
simple $C$-module in the mouth of the tube of rank $p_t$ and
consider $S$ as an object in $\coh\,(\x)$. The category $S^\perp$
right perpendicular to the object $S$ is the full subcategory of
$\coh\,(\x)$ consisting of all $\mathcal{F}\in \coh\,(\x)$
satisfying
$$\Hom_{\x}(S,\mathcal{F})=0=\Ext^1_{\x}(S,\mathcal{F}).$$

By \cite{20}, $S^\perp =\coh\,(\x')$ where $\x' =\x(p',\lambda)$ is
a weighted projective line with weight sequence
$p'=(p_1,p_2,\ldots,p_{t-1},p_t-1)$. Moreover, if $0\to \tau S\to
U\to S\to 0$ is the almost split sequence in $\coh\,(\x)$, then $U$
is a simple object in $S^\perp$ of $\tau'$-period $p_t-1$, where
$\tau' =\tau_{\rD^b(\coh\,(\x'))}$.

\begin{propo}
Let $C=C(p,\lambda)$ be a canonical algebra of tubular type
$p=(p_1,\ldots,p_t)$ and $A=C[P]$ be an extended canonical algebra.
Then $A$ is derived canonical of type
$\bar{p}=(p_1,\ldots,p_{t-1},p_t+1)$.
\end{propo}

\begin{proof}
By (\ref{sect:2.4}), we may choose $P$ to be the simple projective
$C$-module. We shall show that $A$ is quasi-tilted of type
$\bar{p}=(p_1,\ldots,p_{t-1},p_t+1)$, see \cite{21}.

Let $\x$ be a weighted projective line with
$\rD^b(\coh\,(\x))=\rD^b(\mod\,C)$ and let $\bar{\x}$ denote a
weighted projective line of type $\bar{p}$ such that $\coh\,(\x)$ is
the perpendicular category $S^\perp$ formed in $\coh\,(\bar{\x})$
for a simple $S$ from the tube of rank $p_t+1$. Let $U$ be the
middle term of the almost split sequence $0\to \tau S\to U\to S\to
0$ in $\coh\,(\bar{\x})$. Then $U$ is a simple in $S^\perp
=\coh\,(\x)$ belonging to the largest tube in $\coh\,(\x)$.

By hypothesis, $\x$ has tubular type. By \cite{19}, there is a
tilting object $\mathcal{T}$ in $\coh_0(\x)\lor \coh_+(\x)[1]$ such
that $\End\,(\mathcal{T})=C$, where $\coh_0(\x)$ (resp.
$\coh_+(\x)$) denotes the full subcategory of $\coh\,(\x)$ formed by
the sheaves of rank $0$ (resp. positive rank). Therefore,
$\mathcal{T}\cup \{S\}\subset \coh_0(\bar{\x})\lor
\coh_+(\bar{\x})[1]$ is a tilting complex for $\bar{\x}$ whose
endomorphism ring is isomorphic to $C[P]=A$.
\end{proof}

\subsection{Negative Euler characteristic: the wild case}
For negative Euler characteristic the derived category of modules
over an extended canonical algebra $C[P]$ relates to the study of
the $\ZZ$-graded surface singularity $R$ associated with $C$ and the
weighted projective line $\x$ associated to $C$. We refer to
\cite{5, 1, 20} for further details. The weighted projective line
$\x=\x(p,\la)$ for a weight sequence $p=(p_1,\ldots,p_t)$ and a
parameter sequence $\la=(\la_3,\ldots,\la_t)$ was introduced in
\cite{5} by means of the algebra
 $$
S=S(p,\la):=k[x_1,\ldots,k_t]/(x_i^{p_i}=x_2^{p_2}-\la_ix_1^{p_1}),\quad
i=3,\ldots,t.
 $$
The algebra $S$ is naturally graded over the abelian group $\Lp$
with generators $\vx_1,\ldots,\vx_t$ and relations
$p_1\vx_1=\cdots=p_t\vx_t=:\vc$ by giving each $x_i$ the degree
$\vx_i$. Here, $\vc$ is called the canonical element of $\Lp$. The
group $\Lp$ is isomorphic to the direct sum of the group $\ZZ$ of
integers and some finite group. A quick way to arrive at the
category $\coh\,(\x)$ of coherent sheaves on $\x$ is by putting
$$
\coh\,(\x)=\frac{\mod^{\Lp}\hy S}{{\mod^{\Lp}_0}\hy S},
$$
where the quotient category on the right is formed in the sense of
\cite{32} and the categories $\mod^{\Lp}\hy S$ (resp.\
${\mod^{\Lp}_0}\hy S$) are the categories of finitely generated
$\Lp$-graded $S$-modules (resp.\ those of finite length).

The element $\vom=(t-2)\vc-\sum_{i=0}^t \vx_i$ from $\Lp$ is called
the \emph{dualizing element}. Its importance comes from the fact
that Serre duality for $\coh\,(\x)$ holds in the form
$D\Ext^1(X,Y)=\Hom(Y,X(\vom))$, where $X\mapsto X(\vom)$ is the
self-equivalence of $\coh\,(\x)$ induced by grading shift $M\mapsto
M(\vom)$, given by $M(\vom)_{\vx}=M_{\vom+\vx}$.

Assume that $\euler{\x}<0$. Then the \emph{graded surface
singularity $R=R(p,\la)$ attached to the weighted projective line}
$\x$ (or the canonical algebra $C$) with data $(p,\la)$ is defined
as $R=\bigoplus_{n=0}^\infty R_n$, where $R_n=S_{n\vom}$. It follows
immediately that $R$ is a finitely generated, i.e.\ affine,
$k$-algebra where each $R_n$ is finite dimensional over $k$ and,
moreover, $R_0=k$ and $R_1=0$. The next theorem illustrates the role
of $R$, and shows in particular that the algebra $R$ keeps all
information on $\x$. For the proofs we refer to \cite{1, 20}.

\begin{theor}
Assume $\euler{\x}<0$. Then the following holds:

$\mathrm{(a)}$ The algebra $R=R(p,\la)$ is a positively $\ZZ$-graded
isolated surface singularity which is graded Gorenstein of
Gorenstein index $-1$.

$\mathrm{(b)}$ There is a natural equivalence $\coh\,(\x)\ra
\mod^\ZZ\hy R/{\mod^\ZZ_0}\hy R$, induced by restricting the grading
from $\Lp$ to $\ZZ=\ZZ\vom$.

$\mathrm{(c)}$ For $k=\CC$, the algebra $R(p,\la)$ is the positively
$\ZZ$-graded algebra of automorphic forms on the (upper) complex
half-plane $\HHplus$ with respect to the action of a Fuchsian group
$G$ of the first kind of signature $(0;p_1,\ldots,p_t)$. \cuadro
\end{theor}

Concerning (a) we note that --- restricting to the case of Krull
dimension two --- the Gorenstein index $d$ can be defined through
the minimal graded injective resolution $0\ra R\ra E^0\ra E^1\ra E^2
\ra0$ of the $R$-module $R$, where the term $E^2$ is the graded
injective hull of $k(d)$ and generally the \emph{grading shift}
$(n)$ is defined by $M(n)_m=M_{n+m}$. In this situation, Serre
duality holds for $\coh(\x)$ in the form
$D\Ext^1(X,Y)=\Hom(Y,X(-d))$, such that the Auslander-Reiten
translation comes from the grading shift $X\mapsto X(-d)$.

Concerning (c) we remark that $G$ is the orbifold fundamental group
of $\x$, having a presentation $\langle \si_1,\ldots,\si_t\,\mid\,
\si_1^{p_1}=\ldots=\si_t^{p_t}=\si_1\cdots\si_t \rangle$ acting by
covering transformations on the (branched) universal cover $\HHplus$
of $\x$. We refer to \cite{23, 24} for the associated rings of
automorphic forms.

\subsection{} For a variety $X$ Orlov investigated in \cite{33} the
triangulated category $\Dsg(X)$ of the singularities of $X$ defined
as the quotient of the bounded derived category $\rD^b(\coh(X))$ of
coherent sheaves modulo the full subcategory of perfect complexes.
If $X$ is affine with coordinate algebra $R$ this category $\Dsg(R)$
is just the quotient $\rD^b(\mod\hy R)/\rD^b(\proj\hy R)$, where
$\proj\hy R$ is the category of finitely generated projective
$R$-modules. In \cite{34} Orlov further introduced a graded variant
$$
\Dsg^\ZZ(R)=\rD^b(\mod^\ZZ\hy R)/\rD^b(\proj^\ZZ\hy R)
$$
called the \emph{triangulated category of the graded singularity
$R$} which will play a central role in this section.

Under the name \emph{stabilized derived category of $R$} the
categories $\Dsg(R)$ were introduced by Buchweitz in \cite{36}. His
results easily extend to the graded case and yields for an $R$ that
is graded Gorenstein an alternative description of $\Dsg^\ZZ(R)$ as
the \emph{stable category of graded maximal Cohen-Macaulay modules}
$\stabMCM^\ZZ\hy R$. More precisely, he showed that the category
$\MCM^\ZZ\hy R$ of maximal graded Cohen-Macaulay $R$-modules is a
Frobenius-category, hence inducing
--- in Keller's terminology \cite{28} --- on the attached stable
category $\stabMCM^\ZZ\hy R$ of graded maximal Cohen-Macaulay
modules modulo projectives,  the structure of an \emph{algebraic}
triangulated category. For a related approach measuring the
complexity of a singularity by a triangulated category we refer to
Krause's account~\cite{37}.

Let $R=\bigoplus_{n\geq0}R_n$, $R_n=S_{n\vom}$, be the positively
$\ZZ$-graded Gorenstein singularity attached to the weighted
projective line $\x$. It follows from \cite[5.6]{1} that $R$ has
Krull dimension two and Gorenstein index $-1$. We fix some notation:
Let $\Mm=\rD^b(\mod^\ZZ\hy R)$ and $\Mmplus=\rD^b(\mod^{\ZZplus}\hy
R)$. Let $\Ppplus$ be the triangulated subcategory of $\Mmplus$
generated by all $R(-n)$, $n\geq0$ and $\Tt$ its left perpendicular
category $\lperp{\Ppplus}$ formed in $\Mmplus$. Denote further by
$\Ssplus$ the triangulated subcategory of $\Mmplus$ generated by all
$k(-n)$, $n\geq0$ and $\Dd$ its right perpendicular category
$\rperp{\Ssplus}$ formed in $\Mmplus$. Finally let
$\Dd(-1)=\rperp{\Ssplus(-1)}$. Then \cite[2.5]{34} implies the
following proposition.

\begin{propo} Assume that $\euler{\x}<0$ and let $R$ be the
positively $\ZZ$-graded singularity attached to $\x$. Then the
following holds:

$\mathrm{(a)}$ The natural functor $\Tt \incl \Mm\stackrel{q}{\ra}
\Dsing(R)$, where $q$ is the quotient functor, is an equivalence of
triangulated categories.

$\mathrm{(b)}$ The $R$-module $k$ is an exceptional object in $\Tt$
with $\lperp{k}=\Dd(-1)$. Moreover, the category $\rD^b(\coh\,(\x))$
is naturally equivalent to $\Dd(-1)$ under the functor $Y\mapsto
(\bR \Gaplus(Y))(-1)$.
\end{propo}

\begin{proof} For the convenience of the reader we sketch the
argument. Using that $R$ is Gorenstein of Gorenstein index -1, and
invoking Gorenstein duality $\bR\Hom^\bullet_R(-,R)$ of $\Mm$ one
sees that $\rperp{T}\subset \rperp{\Dd(-1)}$ and hence $\Dd(-1)$ is
a full subcategory of $\Tt$. Further we see that
$\lperp{k}=\Dd(-1)$. It is well-known that
$$
\Gaplus:\coh(\x)\ra \mod^{\ZZplus}\hy R,\quad Y\mapsto
\bigoplus_{n=0}^\infty \Hom(\Oo,Y(n))
$$
is a full embedding having sheafification, that is, the quotient
functor $q_+:\mod^{\ZZplus}\hy R\ra \coh(\x)$ as an exact left
adjoint and such that composition $q\,\Gaplus$ is the identity
functor on $\coh(\x)$, compare~\cite[1.8]{5}, \cite[5.7]{1}. It
follows that $\bR\Gaplus:\rD^b(\coh(\x))\ra \Mmplus$ is a full
embedding having $q_+:\Mmplus\ra\Mmplus/\Ssplus$ as a left adjoint,
and $q_+\,\bR\Gaplus=1$.

Since $R$ is positively graded with $R_0=k$, it follows that $k$ is
exceptional in $\Mm$ and hence in $\Mmplus$. Invoking the minimal
graded injective resolution $0\ra R \ra E^0\ra E^1 \ra E^2\ra 0$,
where $E^0$ and $E^1$ are socle-free and $E^2$ is the graded
injective hull of $k(-1)$, it follows that $k$ belongs to $\Tt$ and
then also to $\Dd(-1)$. It is straightforward to check that
$\rperp{\Dd(-1)}$ equals the triangulated subcategory $\langle
k\rangle$ generated by $k$, and hence $\lperp{k}=\Dd(-1)$ in $\Tt$.
\end{proof}

\subsection{Proof of Theorem~\ref{thm:1}} We are now in a position to clarify the structure of the
category $\Dsing^\ZZ\hy (R)$. The result was first observed by
K.~Saito and A.~Takahashi (personal communication); it is not yet
published, and uses the technique of matrix factorizations as in
\cite{35}.

By Proposition~\ref{propo:special} it suffices to show that the left
adjoint $\ell:\Tt\ra \lperp{k}$ to the inclusion $j:\lperp{k}\incl
\Tt$ maps $k$ to a line bundle in $\Dd(-1)=\rD^b(\coh(\x))$ up to
translation in $\Dd(-1)$. We put $A=(\bR\Gaplus(\Oo(-\vom)))(-1)$
and construct a morphism $\ga:k\ra A[1]$ such that
$\Hom(\ga,Y):\Hom(A[1],Y)\ra\Hom(k,Y)$ is an isomorphism for each
$Y\in\lperp{k}$ such that $\ell(k)=A[1]$.

The claim is proved in two steps. Put $R_+=\bigoplus_{n\geq 1}R_n$,
then the exact sequence $0\ra R_+\ra R \ra k$ yields an exact
triangle $R\ra k \stackrel{\al}{\ra} R_+[1]$ in $\Mm$, where
$\Hom(\al,Y)$ is an isomorphism for each $Y\in\lperp{k}$. Note for
this that $R$ belongs to $\lperp{\Dd(-1)}$.

For the next step it is useful to identify the derived category
$\Mmplus$ with the full subcategory of $\rD^b(\Mod^\ZZplus\hy R)$
consisting of all complexes with cohomology in $\mod^{\ZZplus}\hy
R$. Here, $\Mod^\ZZplus\hy R$ denotes the category of \emph{all}
graded $R$-modules. Let $0\ra R(-1)\ra E^0 \ra E^1 \ra E^2 \ra 0$ be
the minimal graded injective resolution of $R(1)$ such that $E^2$
equals the graded injective envelope of $k$. (This uses that $R$ has
Gorenstein index $-1$.) Sheafification yields the minimal injective
resolution $0\ra \Oo(\vom)\ra \tilde{E}^0\ra \tilde{E}^1 \ra0$ of
$\Oo(\vom)$. Accordingly $\bR\Gaplus(\Oo(\vom))$ is given by the
complex
$$
A: \quad \cdots\ra 0 \ra E_+^0(-1) \ra E_+^2(-1) \ra 0 \cdots,
\textrm{ where } E_+^i=\bigoplus_{n\geq0}E_n^i,
$$
whose cohomology is concentrated in degrees zero and one and given
by
$$
\mathrm{H}^0(A)=R_+,\quad \mathrm{H}^1(A)=k(-1).
$$
It follows the existence of an exact triangle
$$
k(-1)[-2]\ra R_+ \stackrel{\be}{\ra} A \ra k(-1)[-1],
$$
in $\Mmplus$ where, by construction, $A$ belongs to $\Dd(-1)$. For
$Y$ from $\Dd(-1)$ we have, in particular, that $Y$ belongs to
$\lperp{k(-1)}$ implying that $\Hom(\be,Y)$ is an isomorphism. To
summarize: The morphism $\ga=[k\stackrel{\al}{\ra} R_+[1]
\stackrel{\be[1]}{\ra} A[1]]$ yields isomorphisms $\Hom(\ga,Y)$ for
each $Y\in \Dd(-1)$. Hence $\ell(k)[-1]=A=\bR\Gaplus(\Oo(\vom))(-1)$
is a line bundle, as claimed.\cuadro

\subsection{The Coxeter-Dynkin algebras of a singularity}
\label{coxeter:dynkin} In the theory of singularities the attached
Coxeter-Dynkin diagrams, see for instance \cite{40,9}, play an
important role, in particular, since they establish a link to Lie
theory.

\begin{defin}
Let $k$ be an algebraically closed field, and $R=R(p,\la)$ be the
$\ZZ$-graded singularity attached to the weighted projective line
$\x(p,\la)$.

(a) By the \emph{Coxeter-Dynkin algebra of hereditary type} we mean
the path algebra $D[p]$ of the hereditary star $[p_1,\ldots,p_t]$
having a unique sink.

(b) By the \emph{Coxeter-Dynkin algebra of canonical type} we mean
the algebra $D(p,\la)$ given in terms of the quiver {\small
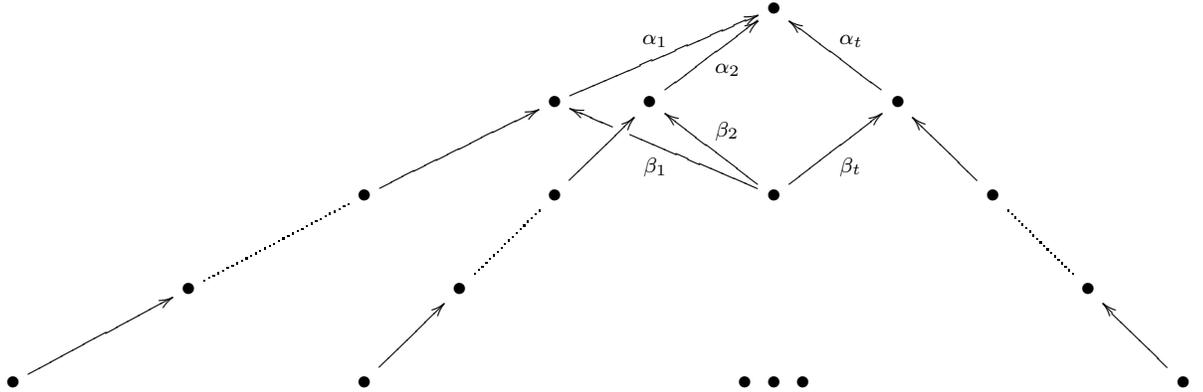
$$\scriptscriptstyle \xymatrix{
 &&&           &              &                &                   &                  &     \bu                                         &                  &          &          &\\
 &&&           &              &                &\bu\ar[urr]^{\al_1}&\bu\ar[ur]_{\al_2}&                                                 &\bu\ar[lu]_{\al_t}&          &          &\\
 &&&           & \bu\ar[urr]   &                &\bu\ar[ur]         &                  &\bu\ar[ur]_{\be_t}\ar[llu]^{\be_1}|!{[u];[lllll]}\hole\ar[lu]_{\be_2}&                  &\bu\ar[lu]&          &\\
 &&\bu\ar@{.}[urr]&           &              &\bu\ar@{.}[ur]  &                   &                  &                                        &                  &          &\bu\ar@{.}[lu]&\\
 \bu\ar[urr]&&&   &\bu\ar[ur]    &                &                   &                  &                                         \bu\;\bu\;\bu        &                  &         &          &\bu\ar[lu]\\
}$$} \noindent with the two relations $\sum_{i=2}^t \al_i\be_i=0$
and $\al_1\be_1=\sum_{i=3}^t \la_i\,\al_i\be_i$.

 (c) By the \emph{Coxeter-Dynkin algebra  of
extended-canonical type} we mean the one-point extension
$\hat{D}(p,\la)$ of the Coxeter-Dynkin algebra $D(p,\la)$ of
canonical type above, introducing a new arrow at the sink vertex and
keeping the relations.
\end{defin}

The link to singularity theory is given by the following well known
result, compare \cite{9,40}.
\begin{theor}
For $k=\C$ the Coxeter-Dynkin diagram of the singularity $R(p,\la)$
is the underlying digraph of the Coxeter-Dynkin algebra

(a) of hereditary type, if $\euler\x>0$, and then $[p_1,\ldots,p_t]$
is Dynkin.

(b) of canonical type, if $\euler\x=0$, and then $p$ is tubular.

(c) of extended canonical type, if $\euler\x<0$.\cuadro
\end{theor}
Recall that the digraph of a finite dimensional algebra $A$ has the
underlying graph of the quiver of $A$ as solid edges and the minimal
number of relations between vertices $i$ and $j$ as dotted edges.

\begin{propo}
Assume the number of weights $>1$ is at least two. Then the
following holds:

(a) The Coxeter-Dynkin algebra $D(p,\la)$ of canonical type is
derived equivalent to the canonical algebra $C(p,\la)$.

(b) The Coxeter-Dynkin algebra $\hat{D}(p,\la)$ of extended
canonical type is derived equivalent to the extended canonical
algebra $A=C[P]$, where $C=C(p,\la)$.
\end{propo}

\begin{proof}
Let $\x$ be the weighted projective line with $t$ weighted points
$x_1,\ldots,x_t$ of weight $p_1,\ldots,p_t$, respectively. For each
$i=1,\ldots,t$ denote by $U_i$ the unique indecomposable sheaf of
length $p_i-1$ concentrated in $x_i$ such there exists an
epimorphism $\Oo\ra U_i$ in $\coh(\x)$, where $\Oo$ is the structure
sheaf on $\x$. Moreover, let
$$
S_i=U_i^{(1)}\subset U_i^{(2)} \subset \cdots \subset
U_i^{(p_i-1)}=U_i
$$
be the complete system of subobjects of $U_i$. It is straightforward
to verify that the object
$$
T=\Oo(\vc) \oplus \bigoplus_{i=1}^t\left(U_i^{(1)}\oplus U_i^{(2)}
\oplus \cdots \oplus U_i^{(p_i-1)}\right)\oplus \Oo(-\vom)[1]
$$
is a tilting object in $\rD^b(\coh(\x))$. Moreover, it is not
difficult to see that the endomorphism ring of $T$ is isomorphic to
the Coxeter-Dynkin algebra $D(p,\la)$.

This proves (a). Assertion (b) follows from (a) noting that
$\Oo(-\vom)$ is a line bundle by applying \cite{13} or arguing as in
Proposition~\ref{propo:special}.
\end{proof}

\subsection{The triangulated category of singularities for
nonnegative Euler characteristic} In this paper we mainly
concentrate on the case $\euler\x<0$. To complete the picture we
review the situation for $\euler\x\geq0$.

Assume that $k=\C$. For $\euler\x>0$ the Dynkin diagram given by the
weight type and the Coxeter-Dynkin diagram of the singularity
$R=R(p,\la)$ agree. Moreover, it is shown in \cite{35} that the
category $\Dsing^\ZZ(R)$ is equivalent to the bounded derived
category $\rD^b(k\Delta)$ of the path algebra of a quiver of the
same Dynkin type.

Next we deal with the case $\euler\x=0$. First, there is no
$\ZZ$-graded Gorenstein algebra $R$ such that
$\mathcal{C}=\mod^\ZZ\hy R/\mod^\ZZ_0\hy R$ is equivalent to
$\coh(\x)$. Assume, indeed, that such an algebra $R$ would exist.
Let $d$ denote its Gorenstein index. Since the Auslander-Reiten
translation $\tau(X)=X(-d)$ has finite order 2, 3, 4 or 6 for any
weighted projective line $\x$ of tubular type, it follows that
$d=0$, hence $\tau$ is the identity, contradiction.

Hence for $\euler\x=0$ it is more natural to investigate the
$\Lp$-graded singularity $S=S(p,\la)$ and its triangulated category
of singularities $\Dsing^{\Lp}(S)$ which in the tubular case is
equivalent to $\rD^b(\coh(\x))$ due to recent work of
Ueda~\cite{38}.

\section{The Coxeter polynomial of an extended canonical algebra} 
\label{sect:cox}
\subsection{} 
\label{sect:4.1} Let $A$ be a finite dimensional $k$-algebra of
finite global dimension. The {\em Coxeter transformation\/}
$\varphi_A\colon K_0(A)\to K_0(A)$ is the automorphism induced by
the Auslander-Reiten translation $\tau_{\rD^b(\mod\,A)}\colon
\rD^b(\mod\hy A)\to \rD^b(\mod\hy A)$. We shall consider $\varphi_A$
as a $n\times n$ integral matrix where $n$ is the rank of $K_0(A)$.
>From the Introduction, we recall that $f_A(T)=\det\,(T\id
-\varphi_A)$ is called the {\em Coxeter polynomial\/} of $A$.

For a one-point extension $A=B[P]$ of an algebra $B$ by an
indecomposable projective a simple calculation shows (see
\cite{3,11}):
$$f_A(T)=(1+T)f_B(T)-Tf_C(T).$$

Of particular interest is the hereditary case where any algebra can
be constructed by repeated one-point extensions using indecomposable
projective modules. In fact, for a star $H$ of type
$[p_1,\ldots,p_t]$ the above formula yields:
$${\textstyle f_{[p_1,\ldots,p_t]}(T):=f_H(T)=
\left[(T+1)-T\sum\limits^t_{i=1}\frac{v_{p_i-1}}{v_{p_i}}\right]
\prod\limits^t_{i=1}v_{p_i}}$$
where we set $v_n(T)=\frac{T^n-1}{T-1}=\sum\limits^{n-1}_{i=0}T^i$.

For the canonical algebra $C=C(p,\lambda)$ of type
$p=(p_1,\ldots,p_t)$, we get
$${\textstyle f_{(p_1,\ldots,p_t)}(T):=
f_C(T)=(T-1)^2\prod\limits^t_{i=1}v_{p_i}(T).}$$
In particular, all the eigenvalues of $\varphi_C$ lie on the unit
circle $\s^1$. For these calculations we refer the reader to \cite{2}.

\begin{lemmas} {\rm \cite{3}.}
An extended canonical algebra $A=C[P]$ where $C=C(p,\lambda)$ with
$p=(p_1,\ldots,p_t)$ has Coxeter polynomial:
$${\textstyle \hat{f}_{(p_1,\ldots,p_t)}(T):=f_A(T)=
(T+1)(T-1)^2\prod\limits^t_{i=1}v_{p_i}(T)-Tf_{[p_1,\ldots,p_t]}(T)}$$
\vglue-26pt \cuadro
\end{lemmas}

\subsection{} 
\label{sect:4.2} For later use we recall some facts on {\em
cyclotomic polynomials}.

The $n$-cyclotomic polynomial $\phi_n(T)$ is inductively defined by
the formula
$${\textstyle T^n-1=\prod\limits_{d\mid n}\phi_d(\mathcal{T}).}$$

Recall that the {\em M\"obius function\/} is defined as follows:
$$\mu (n)=\begin{cases}
0 &\hbox{if $n$ is divisible by a square}\\
(-1)^r &\hbox{if $n=p_1,\ldots p_r$ is a factorization into distinct
primes.}\end{cases}$$ A more explicit expression for the cyclotomic
polynomials is given by:


\begin{lemma}
For each $n\ge 2$, we have
$${\textstyle \phi_n(T)=\prod\limits_{1\le d<n\atop d\mid
n}v_{n/d}(T)^{\mu (d)}}$$
\end{lemma}

\subsection{} 
\label{sect:4.3} Following \cite{3} we say that a polynomial
$p(T)\in\Z [T]$ is {\em represented\/} by $q(T)\in \Z [T]$ if
$$p(T^2)=q^*(T):=T^{\degp\,q}q(T+T^{-1}).$$
The interest in the representability of polynomials is due to the
relation between the set of roots of $p(T)$ and $q(T)$ whenever
$p(T^2)=q^*(T)$. Indeed, in that case $\Root\,p(T)\subset \s^1$
(resp. $\s^1\setminus \{1\}$) if and only if $\Root\,q(T)\subset
[-2,2]$ (resp. $(-2,2)$). In \cite{3} is shown that the Coxeter
polynomial $\hat{f}_{(p_1,\ldots,p_t)}(T)$ of an extended canonical
algebra of type $(p_1,\ldots,p_t)$ is represented by
$${\textstyle q_{(p_1,\ldots,p_t)}(T)=
T(T^2-1)\prod\limits^t_{i=1}v_{p_i}(T)
-\chi_{[p_1,\ldots,p_t]}(T),}$$ where $\chi_{[p_1,\ldots,p_t]}(T)$
is the characteristic polynomial of the {\em adjacency matrix\/} of
the star graph of type $[p_1,\ldots,p_t]$.

Using the above expressions the following was recently shown by the
authors:

\begin{theors} {\rm \cite{3}}.
Let $A\!=\!C[P]$ be an extended canonical algebra of type
$(p_1,\ldots,p_{t-1},p_t+1)$ and $A'$ be an extended canonical
algebra of type $(p_1,\ldots,p_t)$. Then the following holds:

\begin{itemize}
\item[{\rm (a)}] $\varphi_A$ accepts at most $4$ eigenvalues outside
$\s^1$

\item[{\rm (b)}] If\/ $\Root\,f_A\subset \s^1$, then also\/
$\Root\,f_{A'}\subset \s^1$
\end{itemize}
\end{theors}

\begin{trivlist}
\item[\hspace{\labelsep}{\noindent\em Sketch of proof:\/}]
The polynomials $q_{(p_1,\ldots,p_t)}(T)$ satisfy a {\em Chebysheff
recursion formula\/} as follows:
$$q_{(p_1,\ldots,p_t+1)}(T)
=Tq_{(p_1,\ldots,p_t)}(T)-q_{(p_1,\ldots,p_t-1)}(T).$$
A version of Sturm's Theorem applies to assure that for any real
interval $[\alpha ,\beta]$, if $q_{(p_1,\ldots,p_t+1)}(T)$ has roots
$\lambda_1\le \cdots \le \lambda_s$ in $[\alpha,\beta]$, then
$q_{(p_1,\ldots,p_t)}(T)$ has roots $\lambda'_1\le \cdots \le
\lambda'_{s-1}$ in $[\alpha,\beta]$ satisfying
$$\lambda_1\le \lambda'_1\le \lambda_2\le \lambda'_2\le \cdots \le
\lambda_{s-1}\le \lambda'_{s-1}\le \lambda_s.$$
(a) and (b) follow easily from these facts.\cuadro
\end{trivlist}

\subsection{} 
\label{sect:4.4} According to (\ref{sect:4.3}), to prove
Theorem~\ref{thm:2} we need to calculate the minimal weight
sequences $p$ such that $\hat{f}_p(T)$ is not contained in $\s^1$.
This is done by systematically computing the roots of Coxeter
polynomials of extended canonical algebras.

Recall that the {\em spectral radius\/} of the Coxeter
transformation $\varphi_A$ is by definition
$\rho_A(\varphi_A)=\max\,\{|\lambda|\colon \lambda\in
\Root\,f_A(T)\}$. In the following list $A$ is an extended canonical
algebra of weight type $p$. The invariant called {\em Dynkin
index\/} is explained in section~\ref{sect:poinc}.
\newpage
\medskip

{\tiny
\begin{center}
\begin{tabular}{c|c|l|c|c}
&Weight &$\!\!$Irreducible factorization of $f_A(T)$
  &$\!\!\rho (\varphi_A)\!\!$ &$\!\!\!$Dynkin$\!\!\!$\\[8pt]
&$(p_1,\ldots,p_t)$ & & &index\\[8pt]
\hline &&&&\\[-3pt]
$\!\!t\!=\!3\!\!$
  &$(2,3,11)$
    &$\!\!1\!+\!T\!-\!T^3\!-\!T^4\!+\!T^6\!+\!T^7\!+\!T^9\!+\!T^{10}\!-\!T^{12}\!-\!T^{13}\!+\!T^{15}\!+\!T^{16}$
    &$\!\!1.1064\!\!$ &$6$\\[8pt]
  &$(2,4,9)$
    &$\!\!\phi_2\phi_5(T^{10}\!-\!T^9\!+\!T^5\!+\!T\!+\!1)$
    &$\!\!1.1329\!\!$ &$4$\\[8pt]
  &$(2,5,8)$
    &$\!\!1\!+\!T\!+\!T^4\!+\!T^5\!+\!T^6\!+\!2T^8\!+\!T^9\!+\!T^{10}\!+\!T^{11}\!+\!T^{14}\!+\!T^{15}$
    &$\!\!1.1574\!\!$ &$4$\\[8pt]
  &$(2,6,7)$
    &$\!\!1\!+\!T\!+\!T^4\!+\!2T^5\!+\!2T^6\!+\!T^7\!+\!T^8\!+\!2T^9\!+\!2T^{10}\!+\!T^{11}\!+\!T^{14}\!+\!T^{15}$
    &$\!\!1.1669\!\!$ &$4$\\[8pt]
  &$(3,3,8)$
    &$\!\!1\!+\!T\!+\!T^2\!+\!T^5\!+\!2T^6\!+\!3T^7\!+\!2T^8\!+\!T^9\!+\!T^{12}\!+\!T^{13}\!+\!T^{14}$
    &$\!\!1.1498\!\!$ &$3$\\[8pt]
  &$(3,4,7)$
    &$\!\!1\!+\!T\!+\!T^2\!+\!T^3\!+\!T^4\!+\!2T^5\!+\!3T^6\!+\!3T^7\!+\!3T^8\!+\!2T^9\!+\!T^{10}\!+\!T^{11}\!+\!T^{12}\!+\!T^{13}\!+\!T^{14}\!\!$
    &$\!\!1.1847\!\!$ &$3$\\[8pt]
  &$(3,5,6)$
    &$\!\!\phi_3(T^{12}\!+\!T^9\!+\!T^8\!+\!T^7\!+\!T^6\!+\!T^5\!+\!T^4\!+\!T^3\!+\!1)$
    &$\!\!1.1966\!\!$ &$3$\\[8pt]
  &$(4,4,6)$
    &$\!\!\phi_2\phi_4(T^{10}-T^9+T^8+T^6+T^4+T^2-T+1)$
    &$\!\!1.2715\!\!$  &$3$\\[8pt]
  &$(4,5,5)$
    &$\!\!\phi_5(T^{10}\!+\!T^7\!+\!T^6\!+\!T^5\!+\!T^4\!+\!T^3\!+\!1)$
    &$\!\!1.2277\!\!$ &$3$\\[8pt]
\hline &&&&\\[-3pt]
$\!\!t\!=\!4\!\!$
  &$(2,2,2,7)$
    &$\!\!\phi^2_2(T^{10}\!+\!T^6\!+\!T^5\!+\!T^4\!+\!1)$
    &$\!\!1.1670\!\!$ &$2$\\[8pt]
  &$(2,2,3,6)$
    &$\!\!\phi^2_2\phi_3(T^8\!-\!T^7\!+\!T^6\!+\!T^4\!+\!T^2\!-\!T\!+\!1)$
    &$\!\!1.2196\!\!$ &$2$\\[8pt]
  &$(2,3,3,4)$
    &$\!\!\phi_2\phi_4(T^8+T^6+T^5+2T^4+T^3+T^2+1)$
    &$\!\!1.2874\!\!$ &$2$\\[8pt]
  &$(3,3,3,4)$
    &$\!\!\phi^2_3(T^8+T^6+2T^5+2T^3+T^2+1)$
    &$\!\!1.3307\!\!$ &$2$\\[8pt]
\hline &&&&\\[-3pt]
$\!\!t\!=\!5\!\!$
  & & & &$2$\\[8pt]
  &$(2,2,2,2,5)$
    &$\!\!\phi^3_2(T^8\!+\!T^6\!+\!T^3\!+\!2T^4\!+\!T^3\!+\!T^2\!+\!1)$
    &$\!\!1.2874\!\!$ &$2$\\[8pt]
  &$(2,2,2,3,4)$
    &$\!\!\phi^3_2(T^8\!+\!2T^6\!+\!T^5\!+\!3T^4\!+\!2T^2\!+\!1)$
    &$\!\!1.3351\!\!$ &$2$\\[8pt]
  &$(2,2,3,3,3)$
    &$\!\!\phi_2\phi^2_3(T^6\!+\!T^4\!+\!2T^3\!+\!T^2\!+\!1)$
    &$\!\!1.3765\!\!$ &$2$\\[8pt]
\hline &&&&\\[-3pt]
$\!\!t\!=\!6\!\!$
  & & & &$2$\\[8pt]
  &$(2,2,2,2,2,3)$
    &$\!\!\phi^4_2(T^6\!+\!2T^4\!+\!T^3\!+\!2T^2\!+\!1)$
    &$\!\!1.3395\!\!$ &$2$\\[8pt]
\hline &&&&\\[-3pt]
$\!\!t\!=\!7\!\!$
  & & & &$2$\\[8pt]
  &$\!\!\!(2,2,2,2,2,2,2)\!\!\!$
    &$\!\!\phi^6_2(T^4\!-\!T^3\!+\!3T^2\!-\!T\!+\!1)$
    &$\!\!1.5392\!\!$ &$2$
\end{tabular}
\vskip12pt
{\normalsize{\bf Table 1.} Critical weight sequences.} 
\end{center}
}

\medskip

\subsection{} 
\label{sect:4.5} To complete the proof of Theorem~\ref{thm:2} we
shall show that any weight sequence $p'<p$ with $p$ in Table~1 has
all its roots on $\s^1$. This is computed in the following Table~2.

As above $A=C[P]$ is an extended canonical algebra of weight type
$p=(p_1,\ldots,p_t)$. In {\em Table\/}~2, the marks $\bulitito$ and
\cuadrito\ refer to the case $k=\C$: those weight sequences marked
by $\bulitito$ or \cuadrito\ correspond to algebras $R(p,\lambda)$
associated to hypersurface singularities, in those cases
$R(p,\lambda)$ is formally $3$-generated. The marks $\bulitito$
correspond to {\em Arnold's $14$ exceptional unimodal
singularities\/} in the theory or singularities of differentiable
maps \cite{25}. Among those weight sequences $p=(p_1,p_2,p_3)$ (that
is $t=3$), Arnold's singularities are exactly those rings of
automorphic forms having three generators \cite{8}.

\newpage

\begin{center}
{\tiny
\begin{tabular}{rc|l|l|c}
&Weight sequence\ \ &\ \ Factorization of $f_A(T)$\ \ &\ \
Poincar\'e
  series\ \ &\ \ Period of $\varphi_A$\ \ \\[3.5pt]
\hline
&&&&\\[-9pt]
$\bulitito$ &$(2,3,7)$ &$\phi_{42}$ &$(6,14,21)$ $(42)$ &$42$\\[3.5pt]
$\bulitito$ &$(2,3,8)$ &$\phi_2\cdot \phi_{10}\cdot \phi_{30}$
  &$(6,8,15)$ $(30)$ &$30$\\[3.5pt]
$\bulitito$ &$(2,3,9)$ &$\phi_3\cdot \phi_{12}\cdot \phi_{24}$
  &$(6,8,9)$ $(24)$ &$24$\\[3.5pt]
&$(2,3,10)$ &$\phi_2\cdot \phi_{16}\cdot \phi_{18}$
  &$(6,8,9,10)$ $(16,18)$ &$72$\\
\multicolumn{5}{c}{\phantom{\tiny .}\punteada\phantom{\tiny .}}\\
$\bulitito$ &$(2,4,5)$ &$\phi_2\cdot \phi_6\cdot \phi_{30}$
  &$(4,10,15)$ $(30)$ &$30$\\[3.5pt]
$\bulitito$ &$(2,4,6)$ &$\phi^2_2\cdot \phi_{22}$
  &$(4,6,11)$ $(22)$ &$22$\\[3.5pt]
$\bulitito$ &$(2,4,7)$ &$\phi_2\cdot \phi_9\cdot \phi_{18}$
  &$(4,6,7)$ $(18)$ &$18$\\[3.5pt]
&$(2,4,8)$ &$\phi^2_2\cdot \phi_4\cdot \phi_{12}\cdot \phi_{14}$
  &$(4,6,7,8)$ $(12,14)$ &$84$\\
\multicolumn{5}{c}{\phantom{\tiny .}\punteada\phantom{\tiny .}}\\
$\bulitito$ &$(2,5,5)$ &$\phi_5\cdot \phi_{20}$
  &$(4,5,10)$ $(20)$ &$20$\\[3.5pt]
$\bulitito$ &$(2,5,6)$ &$\phi_2\cdot \phi_8\cdot \phi_{16}$
  &$(4,5,6)$ $(16)$ &$16$\\[3.5pt]
&$(2,5,7)$ &$\phi_{11}\cdot \phi_{12}$
  &$(4,5,6,7)$ $(11,12)$ &$132$\\[3.5pt]
&$(2,6,6)$ &$\phi^2_2\cdot \phi_3\cdot \phi_6\cdot \phi_{10}\cdot \phi_{12}$
  &$(4,5,6,6)$ $(10,12)$ &$60$\\
\multicolumn{5}{c}{\phantom{\tiny .}\punteada\phantom{\tiny .}}\\
$\bulitito$ &$(3,3,4)$ &$\phi_3\cdot \phi_{24}$
  &$(3,8,12)$ $(24)$ &$24$\\[3.5pt]
$\bulitito$ &$(3,3,5)$ &$\phi_2\cdot\phi_3\cdot\phi_6\cdot\phi_{18}$
  &$(3,5,9)$ $(18)$ &$18$\\[3.5pt]
$\bulitito$ &$(3,3,6)$ &$\phi^2_3\cdot\phi_{15}$
  &$(3,5,6)$ $(15)$ &$15$\\[3.5pt]
&$(3,3,7)$ &$\phi_2\cdot\phi_3\cdot\phi_4\cdot\phi_{10}\cdot\phi_{12}$
  &$(3,5,6,7)$ $(10,12)$ &$60$\\[3.5pt]
$\bulitito$ &$(3,4,4)$ &$\phi_2\cdot\phi_4\cdot\phi_{16}$
  &$(3,4,8)$ $(16)$ &$16$\\[3.5pt]
$\bulitito$ &$(3,4,5)$ &$\phi_{13}$ &$(3,4,5)$ $(13)$ &$13$\\[3.5pt]
&$(3,4,6)$ &$\phi_2\cdot\phi_3\cdot\phi_9\cdot\phi_{10}$
  &$(3,4,5,6)$ $(9,10)$ &$90$\\[3.5pt]
&$(3,5,5)$ &$\phi_2\cdot\phi_5\cdot\phi_8\cdot\phi_{10}$
  &$(3,4,5,5)$ $(8,10)$ &$40$\\[3.5pt]
$\bulitito$ &$(4,4,4)$ &$\phi^2_2\cdot\phi^2_4\cdot\phi_6\cdot\phi_{12}$
  &$(3,4,4)$ $(12)$ &$12$\\[3.5pt]
&$(4,4,5)$ &$\phi_2\cdot\phi_4\cdot\phi_8\cdot\phi_9$
  &$(3,4,4,5)$ $(8,9)$ &$72$\\
\multicolumn{5}{c}{\phantom{\tiny .}\punteada\phantom{\tiny .}}\\
$\cuadrito$ &$(2,2,2,3)$ &$\phi^2_2\cdot \phi_{18}$
  &$(2,6,9)$ $(18)$ &$18$\\[3.5pt]
$\cuadrito$ &$(2,2,2,4)$ &$\phi^2_2\cdot \phi_{14}$
  &$(2,4,7)$ $(14)$ &$14$\\[3.5pt]
$\cuadrito$ &$(2,2,2,5)$ &$\phi^2_2\cdot\phi_3\cdot\phi_6\cdot \phi_{12}$
  &$(2,4,5)$ $(12)$ &$12$\\[3.5pt]
&$(2,2,2,6)$ &$\phi^2_2\cdot \phi_8\cdot \phi_{10}$
  &$(2,4,5,6)$ $(8,10)$ &$40$\\[3.5pt]
$\cuadrito$ &$(2,2,3,3)$ &$\phi_2\cdot \phi_3\cdot\phi_4\cdot \phi_{12}$
  &$(2,3,6)$ $(12)$ &$12$\\[3.5pt]
$\cuadrito$ &$(2,2,3,4)$ &$\phi^2_2\cdot\phi_5\cdot \phi_{10}$
  &$(2,3,4)$ $(10)$ &$10$\\[3.5pt]
&$(2,2,3,5)$ &$\phi_2\cdot \phi_7\cdot \phi_8$
  &$(2,3,4,5)$ $(7,8)$ &$56$\\[3.5pt]
&$(2,2,4,4)$ &$\phi^2_2\cdot \phi_4\cdot \phi_6\cdot \phi_8$
  &$(2,3,4,4)$ $(6,8)$ &$24$\\[3.5pt]
$\cuadrito$ &$(2,3,3,3)$ &$\phi^2_3\cdot \phi_9$
  &$(2,3,3)$ $(9)$ &$9$\\[3.5pt]
&$(2,3,3,4)$ &$\phi_2\cdot \phi_3\cdot \phi_6\cdot\phi_7$
  &$(2,3,3,4)$ $(6,7)$ &$42$\\[3.5pt]
&$\xy (0,0) *+[F]{(3,3,3,3)}\endxy$ &$\phi_2\cdot\phi^3_3\cdot\phi^2_6$
  &$(2,3,3,3)$ $(6,6)$ &$\infty$\\[3.5pt]
$\cuadrito$ &$(2,2,2,2,2)$ &$\phi^4_2\cdot \phi_{10}$
  &$(2,2,5)$ $(10)$ &$10$\\[3.5pt]
$\cuadrito$ &$(2,2,2,2,3)$ &$\phi^3_2\cdot \phi_4\cdot \phi_8$
  &$(2,2,3)$ $(8)$ &$8$\\[3.5pt]
&$\xy (0,0) *+[F]{(2,2,2,2,4)}\endxy$ &$\phi^2_2\cdot \phi_3\cdot\phi^2_6$
  &$(2,2,3,4)$ $(6,6)$ &$\infty$\\[3.5pt]
&$(2,2,2,3,3)$ &$\phi^2_2\cdot \phi_3\cdot\phi_5\cdot\phi_6$
  &$(2,2,3,3)$ $(5,6)$ &$30$\\[3.5pt]
&$(2,2,2,2,2,2)$ &$\phi^5_2\cdot \phi_4\cdot\phi_6$
  &$(2,2,2,3)$ $(4,6)$ &$12$
\end{tabular}
}
\vskip12pt
{\bf Table 2.} Weights $p$ with $\rho (\varphi_A)=1$.
\end{center}

\medskip

\subsection{Proof of Theorem~\ref{thm:2}} 
(a) follows from (\ref{sect:4.4}) and (\ref{sect:4.5}) using Theorem
(\ref{sect:4.3}). Part (b) is shown in \cite{3}, see
(\ref{sect:4.3}).\cuadro

\section{The Poincar\'e series of an extended canonical algebra} 
\label{sect:poinc}
\subsection{} 
\label{sect:5.1} Let $C=C(p,\lambda)$ be a  wild canonical algebra
with weight sequence $p=(p_1,\ldots,p_t)$. Let $P$ be an
indecomposable projective $C$-module. We define the {\em Poincar\'e
series\/} $\hat{P}_C=\hat{P}(p_1,\ldots,p_t)\in \Z [[T]]$ by
$${\textstyle\hat{P}_C(T)=\sum\limits^\infty_{n=0}\langle
[P],\varphi^n_C[P]\rangle_CT^n.}$$ Recall that
$\varphi^n_C[P]=[\tau^n_{\rD^b(\mod\hy C)}P]$ in $K_0(\rD^b(\mod\hy
C))=K_0(C)$. Moreover, observe that $T+\hat{P}_C(T)$ is the
Hilbert-Poincar\'e series $P_C(T)$, as defined in the Introduction,
for the graded algebra
$${\textstyle
R(p,\lambda)=\bigoplus\limits^\infty_{n=0}\Hom\,(L,\tau^n_{\x}L)}$$
where $\x =\x\,(p,\lambda)$ is a weighted projective line,
$\tau_{\x}$ is the Auslander-Reiten translation in $\coh\,(\x)$ and
$L$ is any rank one bundle, see \cite{2} and \cite{22}. In
particular, $\hat{P}_C(T)$ does not depend on the choice of $P$.

\begin{propo}
With the notation above, let $A\!=\!C[P]$ be an extended canonical
algebra. Then

\begin{itemize}
\item[{\rm (a)}] {\rm \cite[Cor~3.6]{22}:} $f_A(T)=P_C(T)f_C(T)$

\medskip

\item[{\rm (b)}] {\rm \cite[Th.~8.6]{2}:}
$P_C(T)=1+T-T\frac{f_{[p_1,\ldots,p_t]}(T)}{f_{(p_1,\ldots,p_t)}(T)}$

\item[{\rm (c)}]{\rm \cite[Prop.~4.3]{2}:}
$P_C(T)=T+\frac{1}{1-T}+(t-2)\frac{T}{(1-T)^2}-\sum_{i=1}^t
\frac{T}{(1-T)(1-T^{p_i})}$.\cuadro
\end{itemize}
\end{propo}

\subsection{} 
\label{sect:5.2} We recall from \cite{2} the following concepts.

\begin{defins} \cite{2}.
Assume $p=(p_1,\ldots,p_t)$ is a weight sequence of wild type. The
{\em Dynkin label\/} of $p$ is the Dynkin diagram of one of the
extended Dynkin graphs $[2,2,2,2]$, $[3,3,3]$ $[2,4,4]$ or $[2,3,6]$
specified as follows:

\begin{itemize}
\item[{\rm (a)}] if $t\ge 4$, then the label is of type $[2,2,2,2]$

\item[{\rm (b)}] if $t=3$, then the label is of type $[a,b,c]$ if
$[a,b,c]\le [p_1,p_2,p_3]$ and $a+b+c$ is minimal.
\end{itemize}

We say that $p$ has {\em Dynkin index\/} $2$, $3$, $4$ or $6$ if its
Dynkin label is $[2,2,2,2]$, $[3,3,3]$, $[2,4,4]$ or $[2,3,6]$
respectively.
\end{defins}

Consider the graded algebra $R=R(p,\lambda)$. In \cite{2} the {\em
support monoid\/} $M(p)$ was introduced as the set of those $n\in \N$
with $R_n\ne 0$. Clearly, $M(p)$ is an additive semigroup in $\N$
generating $\Z$ as a group.

\begin{propos} {\rm \cite{2}.}
The support monoid $M(p)$ is finitely generated with at most $6$
generators. The smallest element in $M(p)$ is the Dynkin index of
$p$.\cuadro
\end{propos}

\subsection{} 
\label{sect:5.3} In a preliminary version of \cite{2}, the authors
displayed the list of all possible support monoids $M(p)$. This list
of $22$ semigroups is essential for the proofs of
Theorems~\ref{thm:3} and \ref{thm:4} and we reproduce it below. For
a weight sequence $p$ the largest integer $n$ such that $n$ does not
belong to $M(p)$ is called the {\em Frobenius number\/} of $M(p)$
and it is denoted by $\alpha (p)$. Of particular interest is the
fact that $p\le q$ implies $M(p)\subset M(q)$.
\newpage

\vskip5pt 

\begin{center}
{\tiny
\begin{tabular}{||c|c|c||}
\hline
&&\\[-9.5pt]
\hline
\ \ Weight type $p$\ \ &\ \ Frobenius number $\alpha (p)$\ \
  &\ \ generators of $\M(p)$\ \ \\
\hline
&&\\[-9.5pt]
\hline
&&\\[-9.5pt]
$\boldsymbol{(2,3,7)}$ &$43$ &$\{6,14,21\}$\\
&&\\[-11pt]
$(2,3,8)$ &$25$ &$\{6,8,15\}$\\
&&\\[-11pt]
$(2,3,9)$ &$19$ &$\{6,8,9\}$\\
&&\\[-11pt]
$(2,3,10)$ &$13$ &$\{6,8,9,10\}$\\
&&\\[-11pt]
$(2,3,11)$ &$13$ &$\{6,8,9,10,11\}$\\
&&\\[-11pt]
$(2,3,12)$ &$13$ &$\{6,8,9,10,11\}$\\
&&\\[-11pt]
$(2,3,13)$ &$7$ &$\{6,8,9,10,11,13\}$\\
&&\\[-11pt]
$(2,3,\infty)$ &$7$ &$\{6,8,9,10,11,13\}$\\
&&\\[-11pt]
\hline
&&\\[-9.5pt]
\hline
&&\\[-9.5pt]
$\boldsymbol{(2,4,5)}$ &$21$ &$\{4,10,15\}$\\
&&\\[-11pt]
$(2,4,6)$ &$13$ &$\{4,6,11\}$\\
&&\\[-11pt]
$(2,4,7)$ &$9$ &$\{4,6,7\}$\\
&&\\[-11pt]
$(2,4,8)$ &$9$ &$\{4,6,7\}$\\
&&\\[-11pt]
$(2,4,9)$ &$5$ &$\{4,6,7,9\}$\\
&&\\[-11pt]
$(2,4,\infty)$ &$5$ &$\{4,6,7,9\}$\\
\hline
&&\\[-9.5pt]
$(2,5,5)$ &$11$ &$\{4,5\}$\\
&&\\[-11pt]
$(2,5,6)$ &$7$ &$\{4,5,6\}$\\
&&\\[-11pt]
$(2,5,7)$ &$3$ &$\{4,5,6,7\}$\\
&&\\[-11pt]
$(2,5,\infty)$ &$3$ &$\{4,5,6,7\}$\\
\hline
&&\\[-9.5pt]
$(2,6,6)$ &$7$ &$\{4,5,6\}$\\
&&\\[-11pt]
$(2,6,7)$ &$3$ &$\{4,5,6,7\}$\\
&&\\[-11pt]
$(2,6,\infty)$ &$3$ &$\{4,5,6,7\}$\\
\hline
$(2,\infty,\infty)$ &$3$ &$\{4,5,6,7\}$\\
\hline
&&\\[-9.5pt]
\hline
&&\\[-9.5pt]
$\boldsymbol{(3,3,4)}$ &$13$ &$\{3,8\}$\\
&&\\[-11pt]
$(3,3,5)$ &$7$ &$\{3,5\}$\\
&&\\[-11pt]
$(3,3,6)$ &$7$ &$\{3,5\}$\\
&&\\[-11pt]
$(3,3,7)$ &$4$ &$\{3,5,7\}$\\
&&\\[-11pt]
$(3,3,\infty)$ &$3$ &$\{3,5,7\}$\\
\hline
&&\\[-9.5pt]
$(3,3,4)$ &$5$ &$\{3,4\}$\\
&&\\[-11pt]
$(3,3,5)$ &$2$ &$\{3,4,5\}$\\
\hline
&&\\[-9.5pt]
$(4,4,4)$ &$5$ &$\{3,4\}$\\
&&\\[-11pt]
$(4,4,5)$ &$2$ &$\{3,4,5\}$\\
&&\\[-11pt]
$(\infty,\infty,\infty)$ &$2$ &$\{3,4,5\}$\\
\hline
&&\\[-9.5pt]
\hline
&&\\[-9.5pt]
$\boldsymbol{(2,2,2,3)}$ &$7$ &$\{2,9\}$\\
&&\\[-11pt]
$(2,2,2,4)$ &$5$ &$\{2,7\}$\\
&&\\[-11pt]
$(2,2,2,5)$ &$3$ &$\{2,5\}$\\
&&\\[-11pt]
$(2,2,2,\infty)$ &$3$ &$\{2,5\}$\\
\hline
&&\\[-9.5pt]
$(2,2,3,3)$ &$1$ &$\{2,3\}$\\
&&\\[-11pt]
$(2,2,3,\infty)$ &$1$ &$\{2,3\}$\\
\hline
&&\\[-9.5pt]
$(\infty,\infty,\infty,\infty)$ &$1$ &$\{2,3\}$\\
\hline
&&\\[-9.5pt]
$\boldsymbol{(2,2,2,2,2)}$ &$3$ &$\{2,5\}$\\
&&\\[-11pt]
$(2,2,2,2,3)$ &$1$ &$\{2,3\}$\\
&&\\[-11pt]
$(2,2,2,2,\infty)$ &$1$ &$\{2,3\}$\\
\hline
&&\\[-9.5pt]
$(2,2,2,2,2,2)$ &$1$ &$\{2,3\}$\\
&&\\[-11pt]
$(\infty,\infty,\ldots,\infty,\infty)$ &$1$ &$\{2,3\}$\\
\hline
\hline
\end{tabular}
}
\vskip3pt
{\bf Table 3.} Semigroups of $(\N,+)$ with the form $M(p)$.
\end{center}

\subsection{} 
\label{sect:5.4} $\!\!$For a given weight sequence
$p=(p_1,\ldots,p_t)$, the Coxeter polynomials
$f_{(p_1,\ldots,p_t)}(T)$ and $\hat{f}_{(p_1,\ldots,p_t)}(T)$ are
readily computed by (\ref{sect:4.1}). In case
$\Root\,\hat{f}_{(p_1,\ldots,p_t)}(T)\subset \s^1$, then using
(3.2), the Poincar\'e series $P_C(T)$ can be written as a rational
function
$${\textstyle P_C(T)=
\frac{\prod\limits^m_{i=1}(1-T^{c_i})}{\prod\limits^n_{j=1}
(1-T^{d_j})}(1-T)^r}$$ for sequences $(d_1,\ldots,d_n)$ and
$(c_1,\ldots,c_m)$ of natural numbers $\ge 2$ and some $r\in \Z$. In
case $p$ is of wild type, then by (\ref{sect:4.1})
$${\textstyle \hat{f}_{(p_1,\ldots,p_t)}(1)=-f_{[p_1,\ldots,p_t]}(1)
=\left[(t-2)-
\sum\limits^t_{i=1}\frac{1}{p_i}\right]\prod\limits^t_{i=1}p_i>0}$$
and by (\ref{sect:5.1}), $P_C(T)$ has a pole of order $2$ at $T=1$,
that is $m-n+r=-2$. Moreover, developing the series $P_C(T)$ we
readily see that $r\ne 0$ implies that the semigroups $M(p)$ is
$\N$, but \cite[Th.~10.4]{2} claims that
$1+\left[\frac{1}{t-2}\right]\le \alpha (p)$, that is, $\alpha
(p)>1$ and therefore $r=0$. We state these considerations in the
following.

\begin{lemma}
Let $p=(p_1,\ldots,p_t)$ be a weight sequence of wild type and
$C=C(p,\lambda)$ be a canonical algebra. Then
$\Root\,\hat{f}_{(p_1,\ldots,p_t)}(T)\subset \s^1$ if and only if
$R(p,\lambda)$ is formally $n$-generated, that is
$$P_C(T)
=\frac{\prod\limits^{n-2}_{i=1}(1-T^{c_i})}
{\prod\limits^n_{j=1}(1-T^{d_j})}$$
for numbers $c_1,\ldots,c_{n-2}$ and $d_1,\ldots,d_n$, all $\ge 2$,
satisfying
$${\textstyle 1+\sum\limits^n_{j=1}d_j=\sum\limits^{n-2}_{i=1}c_i.}$$
\end{lemma}

\begin{proof}
If $\Root\,\hat{f}_{(p_1,\ldots,p_t)}(T)\subset \s^1$, we showed that
$P_C(T)$ has the desired form. Moreover,
$${\textstyle \sum\limits^n_{j=1}d_j+\deg\,\hat{f}_{(p_1,\ldots,p_t)}(T)=
\sum\limits^{n-2}_{i=1}c_i+ \deg\,f_C(T). }$$ The converse follows
from $\hat{f}_{(p_1,\ldots,p_t)}(T)=P_C(T)f_C(T)$.
\end{proof}

In {\em Table\/}~2, for a weight sequence $p=(p_1,\ldots,p_t)$ of
wild type and $C=C(p,\lambda)$ the corresponding canonical algebra
we have calculated (under the column `Poincar\'e series') the
sequences $(d_1,\ldots,d_n)$ and $(c_1,\ldots,c_{n-2})$
corresponding to
$P_C(T)=\frac{\hat{f}_{(p_1,\ldots,p_t)}(T)}{f_{(p_1,\ldots,p_t)}(T)}$.

\subsection{Proof of Theorem~\ref{thm:3}} 
\label{sect:5.5}
 Let $C=C(p,\lambda)$ be a canonical algebra of wild
type and $A=C[P]$ be an extended canonical algebra. Implications (a)
$\Rightarrow$ (b) $\Rightarrow$ (c) are clear.

(c) $\Rightarrow$ (a): Assume $\Root\,f_A(T)\subset \s^1$. By (5.4),
$P_C(T)$ is $n$-generated. Hence $M(p)$ is at most $n$-generated. By
Table~3, if $n\geq5$ then $p=(2,3,p_3)$ with $p_3\geq11$. But a
calculation shows (see Table~2) that $\mathrm{Root}\,
\hat{f}_{(2,3,11)}(T)$ is not contained in $\s^1$. Then
Theorem~\ref{thm:2} implies that
$\mathrm{Root}\,\hat{f}_{(2,3,p_3)}(T)$ is not contained in $\s^1$
for $p_3\geq11$, which shows that $n\leq4$.


For $k=\C$, the list entries in Table~2 correspond to Fuchsian
singularities which are minimal elliptic as classified in \cite{14}.
These rings are graded complete intersection domains.\cuadro

\subsection{} 
\label{sect:5.6} As a consequence of Theorem~\ref{thm:3} and the
classification given in Table~2 we get the following.

\begin{coro}
Let $C=C(p,\lambda)$ be a canonical algebra of wild type. Let
$A=C[P]$ be an extended canonical algebra. The following are
equivalent:

\begin{itemize}
\item[{\rm (a)}] $\varphi_A$ is periodic.

\item[{\rm (b)}] $\Spec\,\varphi_A\subset \s^1$ and $p$ is not
$(3,3,3,3)$ or $(2,2,2,2,4)$.

\item[{\rm (c)}] $f_A(T)=\prod\limits^m_{i=1}\phi_{s_i}(T)^{e_i}$ for
$1\le s_1<s_2<\cdots <s_m$ and $e_i\ge 1$ ($1\le i\le m$), with
$e_m=1$.\cuadro
\end{itemize}
\end{coro}

\section{Graded integral domains with $3$ homogeneous generators} 

\subsection{} 
The following simple remark is well-known.

\begin{lemma}
Let $R$ be a graded complete intersection integral $k$-algebra of
Krull dimension $2$. Then $R$ is generated by $3$ homogeneous
elements if and only if $R=k[x_1,x_2,x_3]/(f)$ with $\deg\,x_i=d_i$
($1\le i\le 3$) and $f$ a homogeneous prime  polynomial. In this
case the Poincar\'e-Hilbert series of $R$ has the form
$${\textstyle \sum\limits^\infty_{n=0}(\dim_kR_n)T^n
=\frac{1-T^c}{(1-T^{d_1})(1-T^{d_2})(1-T^{d_3})}}$$
for some natural numbers $c,d_1,d_2,d_3$ satisfying
$1+d_1+d_2+d_3=c$.
\end{lemma}

\begin{proof}
Assume we have a graded surjection $g\colon k[x_1,x_2,x_3]\to R$
such that $y_i=g(x_i)$ is homogeneous of degree $d_i$, $1\le i\le
3$. Since $R$ is graded integral, then $I=\ker\,g$ is a prime ideal.
Since $R$ has Krull-dimension two, the ideal $I$ has height one,
hence it is principal. Let $I=(f)$ and $\deg\,(f)=c$. Then the
Poincar\'e series has the desired form and $1+d_1+d_2+d_3=c$ by
(\ref{sect:4.4}).
\end{proof}

\subsection{Proof of Theorem~\ref{thm:4}} 
Let $C=C(p,\lambda)$ be a wild canonical algebra and $A=C[P]$ a
corresponding extended canonical algebra.

(a) $\Rightarrow$ (b): Assume $R(p,\lambda)$ is formally
$3$-generated. By (\ref{sect:5.4}), $\Root\,f_A\subset \s^1$. The
result follows from the list given in Table~2.

(b) $\Rightarrow$ (a): follows as above from Table~2 and
(\ref{sect:5.4}).

Assume $R=R(p,\lambda)$ is formally $3$-generated with
$${\textstyle \sum\limits^\infty_{n=0}(\dim_kR_n)T^n=
\frac{1-T^c}{(1-T^{d_1})(1-T^{d_2})(1-T^{d_3})}}$$
with $(d_1,d_2,d_3),(c)$ according to Table~2. We shall consider two
distinguished situations:

if $t=3$, then $R$ is a quasi-homogeneous complete intersection of
the form\break $k[x_1,x_2,x_3]/(f)$ with $\deg\,x_i=d_i$ and $f$ a
homogeneous relation as displayed in (\ref{sect:6.3}). The case
$t\ge 4$ and $k=\C$ is considered in (\ref{sect:6.4}).\cuadro

\subsection{Theorem} {\rm \cite{1}.}\label{sect:6.3} 
{\it Let $C=C(p,\lambda)$ be a canonical algebra of wild weight type
$p=(p_1,p_2,p_3)$ such that the graded algebra $R=R(p,\lambda)$ is
formally $3$-generated. Then $R$ has the form
$$R=k[x,y,z]=k[X,Y,Z]/(F)$$
where the relation $F$, the degree triple $\deg\,(x,y,z)$ and
$\deg\,(F)$ are displayed in Table~4:
}

\medskip

\begin{center}
\begin{tabular}{r|c|c|c|c|c|}
&$p$ &$\deg\,(x,y,z)$ &relation $F$ &$\deg\,(F)$ &Name\\
\hline
&&&&&\\[-12pt]
Index $=6$
&$(2,3,7)$ &$(6,14,21)$ &$Z^2+Y^3+X^7$ &$42$ &$E_{12}$\\
&$(2,3,8)$ &$(6,8,15)$ &$Z^2+X^5+XY^3$ &$30$ &$Z_{11}$\\
&$(2,3,9)$ &$(6,8,9)$ &$Z^2+XZ^2+X^4$ &$36$ &$Q_{10}$\\
\hline
&&&&&\\[-12pt]
\hline
$4$
&$(2,4,5)$ &$(4,10,15)$ &$Z^2+Y^3+X^5Y$ &$30$ &$E_{13}$\\
&$(2,4,6)$ &$(4,6,11)$ &$Z^2+X^4Y+ZY^3$ &$22$ &$Z_{12}$\\
&$(2,4,7)$ &$(4,6,7)$ &$Y^3+X^3Y+XZ^2$ &$18$ &$Q_{11}$\\
&$(2,5,5)$ &$(4,5,10)$ &$Z^2+Y^2Z+X^5$ &$20$ &$W_{12}$\\
&$(2,5,6)$ &$(4,5,6)$ &$XZ^2+Y^2Z+X^4$ &$16$ &$S_{11}$\\
\hline
&&&&&\\[-12pt]
\hline
$3$
&$(3,3,4)$ &$(3,8,12)$ &$Z^2+Y^3+X^4Z$ &$24$ &$E_{14}$\\
&$(3,3,5)$ &$(3,5,9)$ &$Z^2+XY^3+X^3Z$ &$18$ &$Z_{13}$\\
&$(3,3,6)$ &$(3,5,6)$ &$Y^3+X^3Z+XZ^2$ &$15$ &$Q_{12}$\\
&$(3,4,4)$ &$(3,4,8)$ &$Z^2-Y^2Z+X^4Y$ &$16$ &$W_{13}$\\
&$(3,4,5)$ &$(3,4,5)$ &$X^3Y+XZ^2+Y^2Z$ &$13$ &$S_{12}$\\
&$(4,4,4)$ &$(3,4,4)$ &$X^4-YZ^2+Y^2Z$ &$12$ &$U_{12}$\\
\hline
\end{tabular}\lower3.6truecm\hbox{$\qed$}
\vskip11pt
{\bf Table 4.}
\end{center}

\medskip

As observed in \cite{1}, these $14$ equations are equivalent to
Arnold's exceptional unimodal singularities. The equations are
slightly different to those of the singularity theory classification,
but equivalent for $k=\C$.

\subsection{} \label{sect:6.4} 
In view of the identification of $R(p,\lambda)$ with a ring of
automorphic forms (\ref{sect:5.5}) in case $k=\C$, we get:

\begin{theors} {\rm \cite{10}.}
Let $C=C(p,\lambda)$ be a canonical algebra of wild type
$p=(p_1,\ldots,p_t)$ with $t\ge 4$ over the complex numbers $\C$ and
$R=R(p,\lambda)$ be the associated graded algebra. Then the following
are equivalent:

\begin{itemize}
\item[{\rm (a)}] $R(p,\lambda)$ is formally $3$-generated;

\item[{\rm (b)}] $9\le \sum\limits^t_{i=1}p_i\le 11$

\item[{\rm (c)}] there is a parameter sequence
$\lambda'=(\lambda'_3,\ldots,\lambda'_t)$ such that the algebra
$R(p,\lambda')$ is of the form
$$k[x,y,z]=k[X,Y,Z]/(F)$$
where the relation $F$, the degree sequence $\deg\,(x,y,z)$ and
$\deg\,(F)$ are displayed below:
\end{itemize}
\end{theors}

\medskip

\begin{center}
\begin{tabular}{r|c|c|c|c|c|}
&$p$ &$\deg\,(x,y,z)$ &relation $F$ &$\deg\,(F)$ &Name\\
\hline
&&&&&\\[-12pt]
$t=4$
&$(2,2,2,3)$ &$(2,6,9)$ &$Z^2+Y^3+X^9$ &$18$ &$J_{3,0}$\\
&$(2,2,2,4)$ &$(2,4,7)$ &$Z^2+XY^3+X^7$ &$14$ &$Z_{1,0}$\\
&$(2,2,2,5)$ &$(2,4,5)$ &$Y^3+XZ^2+X^6$ &$12$ &$Q_{2,0}$\\
&$(2,2,3,3)$ &$(2,3,6)$ &$Z^2+Y^4+X^6$ &$12$ &$W_{1,0}$\\
&$(2,2,3,4)$ &$(2,3,4)$ &$Y^2Z+XZ^2+X^5$ &$10$ &$S_{1,0}$\\
&$(2,3,3,3)$ &$(2,3,3)$ &$Z^3+Y^3+X^3Y$ &$9$ &$U_{1,0}$\\
\hline
&&&&&\\[-12pt]
\hline
$t=5$
&$(2,2,2,2,2)$ &$(2,2,5)$ &$Z^2+Y^5+X^5$ &$10$ &$NA^1_{0,0}$\\
&$(2,2,2,2,3)$ &$(2,2,3)$ &$YZ^2+Y^4+X^4$ &$8$ &$VNA^1_{0,0}$\\
\hline
\end{tabular}\lower2.07truecm\hbox{$\qed$}
\vskip12pt
{\bf Table 5.}
\end{center}

\bigskip

\bigskip

\direc{H. Lenzing\cr Institut f\"ur Mathematik\cr Universit\"at
Paderborn\cr 33095 Paderborn\cr Germany\cr
helmut@math.uni-paderborn.de\cr \hbox{\quad}\cr
\hbox{\quad}}{J.$\!$~A. de la Pe\~na\cr Instituto de
Matem\'aticas\cr Universidad Nacional Aut\'onoma de M\'exico\cr
Ciudad Universitaria\cr M\'exico 04510, D.$\!$ F.\cr M\'exico\cr
jap@matem.unam.mx\cr}

\end{document}